\documentclass[11pt]{article}
\usepackage[margin=1in,nohead]{geometry}
\usepackage{latexsym}
\usepackage{epsfig,graphicx}
\usepackage{amsmath,amsthm,amssymb}
\usepackage{color}

\def\tick{\text{Yes}}

\begin{document}


\newcommand{\bfm}[1]{\mbox{\boldmath $#1$}}
\newcommand{\reals}{\mbox{\bfm{R}}}
\def\grad{\bigtriangledown}
\def\Hes{\hbox{Hes}}
\def\for{\hbox{for }}
\def\Prob{\hbox{Pr}}
\newcounter{rot}
\newcommand{\bignote}[1]{\vspace{0.25 in}\fbox{\parbox{6in}{#1}}%
\marginpar{\fbox{{\bf {\large Q \therot}}}}\addtocounter{rot}{1}}
\newcommand{\mnote}[1]{\marginpar{\footnotesize\raggedright\linepenalty=5000 #1}}
\def\bB{{\bf C}}
\def\cL{{\cal L}}
\def\bh{{\bf h}}
\def\bu{{\bf u}}
\def\bv{{\bf v}}
\def\bg{{\bf g}}
\def\bx{{\bf x}}
\def\by{{\bf y}}
\def\bz{{\bf z}}
\def\cT{{\cal T}}
\def\FP{FIND}
\def\prob{\hbox{Pr}}
\def\Min{\hbox{Min}}
\def\a{\alpha} \def\b{\beta} \def\d{\delta} \def\D{\Delta}
\def\e{\epsilon} \def\f{\phi} \def\F{{\Phi}} \def\g{\gamma}
\def\G{\Gamma} \def\k{\kappa}\def\K{\Kappa}
\def\bw{{\bf w}}
\def\z{\zeta} \def\th{\theta} \def\TH{\Theta}  \def\l{\lambda}
\def\La{\Lambda} \def\m{\mu} \def\n{\nu} \def\p{\pi}
\def\r{\rho} \def\R{\Rho} \def\s{\sigma} \def\S{\Sigma}
\def\t{\tau} \def\om{\omega} \def\OM{\Omega}
\newcommand{\bP}[1]{{\bf P#1}}

\newtheorem{lemma}{Lemma}
\newtheorem{theorem}{Theorem}
\newtheorem{corollary}[theorem]{Corollary}
\newtheorem{Remark}{Remark}
\newtheorem{definition}{Definition}
\newcommand{\proofstart}{{\bf Proof\hspace{2em}}}
\newcommand{\proofend}{\hspace*{\fill}\mbox{$\Box$}}
\newtheorem{remark}{Remark}
\newtheorem{proposition}{Proposition}
\newtheorem{claim}{Claim}
\newcommand{\ooi}{(1+o(1))}
\newcommand{\ul}[1]{\mbox{\boldmath$#1$}}
\newcommand{\ol}[1]{\overline{#1}}
\newcommand{\wh}[1]{\widehat{#1}}
\newcommand{\wt}[1]{\widetilde{#1}}

\newcommand{\rdup}[1]{{\lceil #1 \rceil }}
\newcommand{\rdown}[1]{{\lfloor #1 \rfloor}}

\newcommand{\brac}[1]{\left(#1\right)}
\newcommand{\sbrac}[1]{\left({\scriptstyle #1}\right)}
\newcommand{\tbrac}[1]{\left({\textstyle #1}\right)}
\newcommand{\smorl}[1]{{\scriptstyle #1}}
\newcommand{\sbfrac}[2]{\left(\frac{\scriptstyle #1}{\scriptstyle #2}\right)}
\newcommand{\sfrac}[2]{\frac{\scriptstyle #1}{\scriptstyle #2}}
\newcommand{\bfrac}[2]{\left(\frac{#1}{#2}\right)}
\def\half{\sfrac{1}{2}}
\def\tO{\tilde{O}}
\newcommand{\rai}{\rightarrow \infty}
\newcommand{\ra}{\rightarrow}

\def\Inf{\hbox{Inf}}
\def\sm{\setminus}
\def\seq{\subseteq}
\def\es{\emptyset}
\newcommand{\ind}[1]{\mbox{{\large 1}} \{#1\}}

\def\E{\mbox{{\bf E}\;}}
\def\Pr{\mbox{{\bf Pr}}}
\def\whp{\mbox{${\bf whp}\;$}}
\def\Max{\hbox{MAX}}
\def\OPT{{\bf OPT}}
\def\tr{\hbox{Tr}}
\def\vol{\mbox{Vol}}
\def\lnorm{\left| \left| }
\def\rnorm{\right| \right|}
\def\var{\mbox{{\bf Var}}}

\begin{center}
{\Large\bf Two new Probability inequalities and Concentration Results}\\
Ravindran Kannan\\
Microsoft Research Labs., India
\end{center}

\section{Introduction}

The study of stochastic combinatorial problems as well as Probabilistic
Analysis of Algorithms are among the many subjects which use concentration inequalities. A central
concentration inequality is the H\"offding-Azuma (H-A) inequality:
For real-valued random variables $X_1,X_2,\ldots X_n$ satisfying
respectively absolute bounds and the Martingale (difference) condition:
$$|X_i|\leq 1\quad ;\quad E(X_i|X_1,X_2,\ldots X_{i-1})=0,$$
the H-A inequality asserts the following tail bound:
$\prob\left( \left| \sum_{i=1}^n X_i\right|\geq t\right) \leq c_1e^{-c_2t^2/n},$
for some constants $c_1,c_2$ (which are the tails of $N(0,n)$, the standard normal density with variance $n$, but for constants.)
Here, we present two theorems both of which considerably weaken the assumption of
an absolute bound, as well as the Martingale condition,
while retaining the strength of the conclusion.
As consequences of our theorems, we derive new concentration results for
many combinatorial problems.

Our Theorem 1 is simply stated. It weakens
the absolute bound of 1 on $|X_i|$ to a weaker condition than
a bound of 1 on some moments (upto the $m$ th moment) of $X_i$.
It weakens the Martingale difference assumption to
requiring that certain correlations be non-positive. The conclusion
upper bounds $E(\sum_{i=1}^n X_i)^m$ (essentially) by the $m$ th moment of $N(0,n)$; it will be easy to get tail bounds from these moment bounds. Note that
both the hypotheses and the conclusion involve bounds on moments upto the same $m$; so finite
moments are sufficient to get some conclusions, unlike in H-A as well as Chernoff bounds
in both of which, one uses the absolute bound to get a bound on $E(e^{X_i})$. Note that
if $X_i$ have power law tails (with only finite moments bounded), no automatic bound on $E(e^{X_i})$ is available.
But, both H-A inequality and Chernoff bounds follow as very special cases of our Theorem 1.

The study of the minimum length of a Hamilton tour through $n$ random points
chosen in i.i.d. trials from the uniform density in the unit square,
was started by the seminal work of Bearwood, Halton and Hammersley \cite{bhh}.
The algorithmic question - of finding an approximately optimal Hamilton tour
in this i.i.d. setting was tackled by Karp \cite{karp2} - and his work not only pioneered the
field of Probabilistic Analysis of Algorithms, but also inspired later TSP algorithms for deterministic inputs, like Arora's \cite{arora}.
Earlier hard concentration results for the minimal length of a Hamilton tour in the i.i.d. case were made easy by Talagrand's inequality \cite{talagrand}.
But all these concentration results for the Hamilton tour problem as well as many other
combinatorial problems \cite{steele}
make crucial use of the fact that the points are i.i.d. and so
random variables like the number of points in a region in the unit square
are very concentrated - have exponential tails. In the modern setting, heavier
tailed distributions are of interest. There are many models of what ``heavy-tailed'' distribution should mean; this is not the subject of this paper.
But as we will see, our theorems are amenable to ``bursts in space'', where each region of space chooses (independently) the number of points that fall in it, but then may choose that many points possibly adversarially; further, the number of points may have power-law tails instead of exponential tails. In other problems, one may have ``bursts in time'', where, each time unit may choose from a power-law tailed distribution the number of arrivals/new items/jobs.

Using Theorem 1, we are able to prove as strong concentration as was known for the i.i.d. case of TSP (but for
constants), but, now allowing bursts in space. We do the
same for the minimum weight spanning tree problem as well.
We then consider random graphs where edge probabilities are not equal.
We show a concentration result for the chromatic number (which has been
well-studied under the traditional model with equal edge probabilities.)
In these cases, we use the traditional Doob Martingale construction to first
cast the problem as a Martingale problem. The moment conditions needed for the hypotheses
of our theorems follow naturally.

But an application where we do not
have a Martingale, but do have the weaker hypothesis of Theorem 1 is
when we pick a random vector(s) of unit length as in the well-known Johnson-Lindenstrauss (JL) Theorem on Random Projections. Using Theorem 1, we prove a more general theorem than JL where heavier-tailed distributions are allowed.

A further weakening of the hypotheses of H-A is
obtained in our Main Theorem - Theorem (\ref{mainthm}) whose proof is more complicated.
In Theorem (\ref{mainthm}), we use information on conditional moments of $X_i$ conditioned
on ``typical values'' of $X_1+X_2+\ldots +X_{i-1}$ as well as the ``worst-case'' values.
This is very useful in many contexts as we show.
Using Theorem 2, we settle the concentration question for (the discrete case of) the
well studied stochastic bin-packing problem \cite{coffman}
proving concentration results which we show are best possible. Here, we prove a bound on the
variance of $X_i$ using Linear Programming duality; we then exploit a feature of Theorem (\ref{mainthm}) (which is also present in Theorem (\ref{mainthm-special})): higher moments have lower weight in our bounds, so for bin-packing, it turns out that higher moments don't need to be carefully bounded. This feature is also used for the next application
which is the well-studied problem of proving concentration for the number $X$ of triangles in the standard random graph $G(n,p)$. While many papers have proved good tail bounds for large deviations, we prove here the first sub-Gaussian tail bounds for all values of $p$ - namely that $X$ has $N(0,\var X)$ tails for deviations upto $(np)^{7/4}$ (see Definition (\ref{normal-tails})). [Such sub-Gaussian bounds were partially known for the easy case when $p\geq 1/\sqrt n$, but not for the harder case of smaller $p$.]
We also give a proof of concentration for the longest increasing subsequence problem.
It is hoped that the theorems here will provide a tool to deal with
heavy-tailed distributions and inhomogeneity in other situations as well.

There have been many sophisticated probability inequalities.
Besides H-A (see McDiarmid \cite{mcdiarmid} for many useful extensions) and Chernoff, Talagrand's inequality already referred to (\cite{talagrand}) has numerous applications. Burkholder's  inequality for Martingales \cite{burkholder} and many later developments
give bounds based on finite moments.
A crucial point here is that unlike the other inequalities, different moments
have different weights in the bounds (the second moment has the highest) and this
helps get better tail bounds. We will discuss comparisons of our results with these
earlier results in section \ref{comparisons}. But one more note is timely here: many previous inequalities have also used Martingale bounds after excluding ``atypical'' cases. But usually, they insist on an absolute bound in the typical case, whereas, here we only insist on moment bounds. It is important to note that many (probably all) individual pieces of our approach have been used before; the contribution here is in carrying out the particular  combination of them which is then able to prove results for a wide range applications.

\section{Theorem 1}

In theorem (\ref{mainthm-special}) below, we weaken the absolute bound $|X_i|\leq 1$ of
H-A to (\ref{Ml-special}). Since this will be usually applied with $n\geq m$, (\ref{Ml-special}) will be weaker than $E(X_i^l|X_1+X_2+\ldots +X_{i-1})\leq 1$ which is in turn weaker than the absolute bound - $|X_i|\leq 1$. We replace the Martingale difference condition $E(X_i|X_1,X_2,\ldots X_{i-1})=0$ by the obviously weaker condition (\ref{snc}) which we will call {\it strong negative correlation}; it is only required for odd $l$ which we see later relates to negative correlation.
Also, we only require these conditions for all $l$ upto a certain even $m$. We prove a bound on the (same) $m$ (which is even) th moment of
$\sum_{i=1}^nX_i$. Thus, the higher the moment bounded by the hypothesis, the higher the moment bounded
by the conclusion. This in particular will allow us to handle things like ``power-law'' tails. The following
definition will be useful to describe tail bounds.

\begin{definition}\label{normal-tails}
Let $a,\sigma$ be positive reals.
We say that a random variable $X$ has
$N(0,\sigma^2)$ tails upto $a$ if there exist constants $c,c'$ such that
for all $t\in [0,a]$, we have
$$\prob \left( |X-EX|\geq t\right)\leq c' e^{-ct^2/\sigma^2}.$$
\end{definition}
Here there is a hidden parameter $n$ (which will be clear from the context) and
the constants $c,c'$ are independent of $n$, whereas $a,\sigma$ could depend on $n$.

\begin{theorem}\label{mainthm-special}
Let $X_1,X_2,\ldots X_n$ be real valued random variables and $m$ an
even positive integer satisfying the following for $i=1,2,\ldots n$\footnote{
$E$ will denote the expectation of the entire expression which follows.}:
\begin{eqnarray}
&EX_i(X_1+X_2+\ldots X_{i-1})^l\leq 0\; , l<m, \hbox{odd}.\label{snc}\\
&E(X_i^l|X_1+X_2+\ldots +X_{i-1})\leq \left( {n\over m}\right)^{(l-2)/2}\; l!\; ,
l\leq m, \hbox{ even}.\label{Ml-special}
\end{eqnarray}
Then, we have
\begin{align*}
& E\left( \sum_{i=1}^n X_i\right)^m  \leq (48nm)^{m/2}.\\
& \sum_{i=1}^n X_i \text{ has } N(0,n) \text{ tails upto }\sqrt{nm}.
\end{align*}
\end{theorem}

\begin{remark} Since under the hypothesis of (H-A), (\ref{snc}) and (\ref{Ml-special})
hold for all $m$, (H-A) follows from the last statement of the theorem. We will also show that Chernoff bounds follow as a simple corollary.
\end{remark}

\begin{remark} Note that for the upper bound in (\ref{Ml-special}), we have $$\left[ \left( {n\over m}\right)^{{l\over 2}-1}l!\right]^{1/l}\approx \left( {n\over m}\right)^{{1\over 2}-{1\over l}}l.$$
The last quantity is an increasing function of $l$ when $n\geq m$, which will hold in most applications.
Thus the requirements on $\left( E(X_i^l|X_1+X_2+\ldots +X_{i-1})\right)^{1/l}$
are the ``strongest'' for $l=2$ and the requirements get progressively ``weaker'' for higher moments. This will be useful, since, in applications, it will be easier to bound the second moment than the higher ones. The same qualitative aspect also holds for the Main Theorem.
\end{remark}

\begin{remark}\label{decoupling}
Here, we give one comparison of Theorem (\ref{mainthm-special}) with perhaps the closest result to it in the literature, namely a result proved by de la Pe\~na ((1.7) of \cite{lapena-2} - slightly specialized to our situation) which asserts: If $X_1,X_2,\ldots X_n$ is a Martingale difference sequence with $E(X_i^2|X_1,X_2,\ldots ,X_{i-1})\leq 2$ for all $i$ and $E(X_i^l|X_1,X_2,\ldots ,X_{i-1})\leq (l!/2)\alpha^{(l/2)-1}$, for {\bf all positive even integers} $l$, where $\alpha$ is some fixed real, then
$$\prob\left( \sum_{i=1}^nX_i\geq t\right) \leq \exp\left( -{ct^2\over n+\sqrt\alpha t}\right).$$
It is easy to see that this implies $N(0,n)$ tails upto $n/\sqrt\alpha$.

Setting $\alpha={n\over m}$, the hypothesis of Theorem (\ref{mainthm-special}) implies \cite{lapena-2}'s hypothesis {\bf upto $l\leq m$}, not for all $l$ as required there. Were we to be given this hypothesis for all $l$ and furthermore assume $X_i$ are Martingale differences (rather than the more general (\ref{snc}) condition), then since $n/\sqrt\alpha=
\sqrt{nm}$, we would get the same conclusion as Theorem (\ref{mainthm-special}). \cite{lapena-2}'s result is stronger in other directions (which we won't discuss here), but, a main point of our theorem is to assume only finite moments since we would like to deal with long-tailed distributions. Further, note that we can apply our theorem with $m=O(\sqrt n)$, whence, \ref{Ml-special} allows moment bounds to grow with $n$ unlike \cite{lapena-2}.
\end{remark}

\proofstart
Let $M_l=\Max_{i=1}^n E(X_i^l|X_1+X_2+\ldots +X_{i-1})$ for even $l\leq m$.
For $1\leq i\leq n$ and $q\in \{ 0,2,4,\ldots m-2,m\}$, define
$$f(i,q) = E\left( \sum_{j=1}^i X_j\right)^q.$$
Using the two assumptions, we derive
the following recursive inequality for $f(n,m)$,
which we will later solve (much as one does in a Dynamic Programming
algorithm):
\begin{equation}\label{evenM-special}
f(n,m) \leq f(n-1,m)+{11\over 5}\sum_{t\in \{ 2,4,6,\ldots m\}}{m^t\over t!}M_tf(n-1,m-t),
\end{equation}
{\bf Proof of (\ref{evenM-special})}:
Let $A=X_1+X_2+\ldots X_{n-1}$. Let $a_l={m^l\over l!}E|X_n|^l|A|^{m-l}$.
Expanding $(A+X_n)^m$, we get
\begin{align}\label{expand}
E(A+X_n)^m&\leq EA^m+mEX_nA^{m-1}+\sum_{l=2}^ma_l.
\end{align}
Now, we note that $EX_nA^{m-1}\leq 0$ by hypothesis (\ref{snc}) and so the
second term may be dropped. [In fact, this would be the only use of the Martingale
difference condition if we had assumed it; we use SNC instead, since it clearly suffices.]
We will next bound the ``odd terms'' in terms of the two
even terms on the two sides using a simple ``log-convexity'' of moments argument.
For odd $l\geq 3$, we have
\begin{align*}
&E|X_n|^l|A|^{m-l}\leq E\left( X_n^{l+1}A^{m-l-1}X_n^{l-1}A^{m-l+1}\right)^{1/2}
\leq (E(X_n^{l+1}A^{m-l-1}))^{1/2}(E(X_n^{l-1}A^{m-l+1}))^{1/2}\\
&\hbox{Also, }{1\over l!}\leq {6\over 5}{1\over\sqrt {(l+1)!}}{1\over\sqrt {(l-1)!}}
\end{align*}
So, $a_l$ is at most 6/5 times the geometric mean of $a_{l+1}$ and $a_{l-1}$ and
hence is at most 6/5 times their arithmetic mean.
Plugging this into (\ref{expand}), we get
\begin{equation}\label{evenM}
E(\sum_{i=1}^n X_i)^m\leq
EA^m+{11\over 5}( a_2+a_4+\ldots +a_m)
\end{equation}
Now, we use the standard trick of ``integrating
over'' $X_n$ first and then over $A$ (which is also crucial for proving H-A) to get for even $l$:
$EX_n^lA^{m-l}=E_A \left( A^{m-l} E_{X_n} (X_n^l|A)\right)\leq M_lEA^{m-l}$
which yields
(\ref{evenM-special}).

We view (\ref{evenM-special}) as a recursive inequality for $f(n,m)$. We will use this
same inequality for the proof of the Main theorem, but there we use an inductive proof; here, instead, we will now ``unravel'' the recursion to solve it. [Note that we cannot use induction since we only
know the upper bound involving $(n/m)^{(l/2)-1}$ on the moments (as in the hypothesis of the
theorem) and as $n$ decreases for induction, this bound gets tighter.]
Note that the dropping the $EX_nA^{m-1}$ ensured that the coefficient
of $EA^m$ is 1 instead of the 11/5 we have in front of the other terms. This is important:
if we had 11/5 instead, since the term does not reduce $m$, but only $n$, we would get a
$(11/5)^n$ when we unwind the recursion. This is no good; we can get $m$ terms in the exponent
in the final result, but not $n$.

Imagine a directed graph (see figure Recursion Tree)
\begin{figure}
\begin{center}
\includegraphics[width=4in]{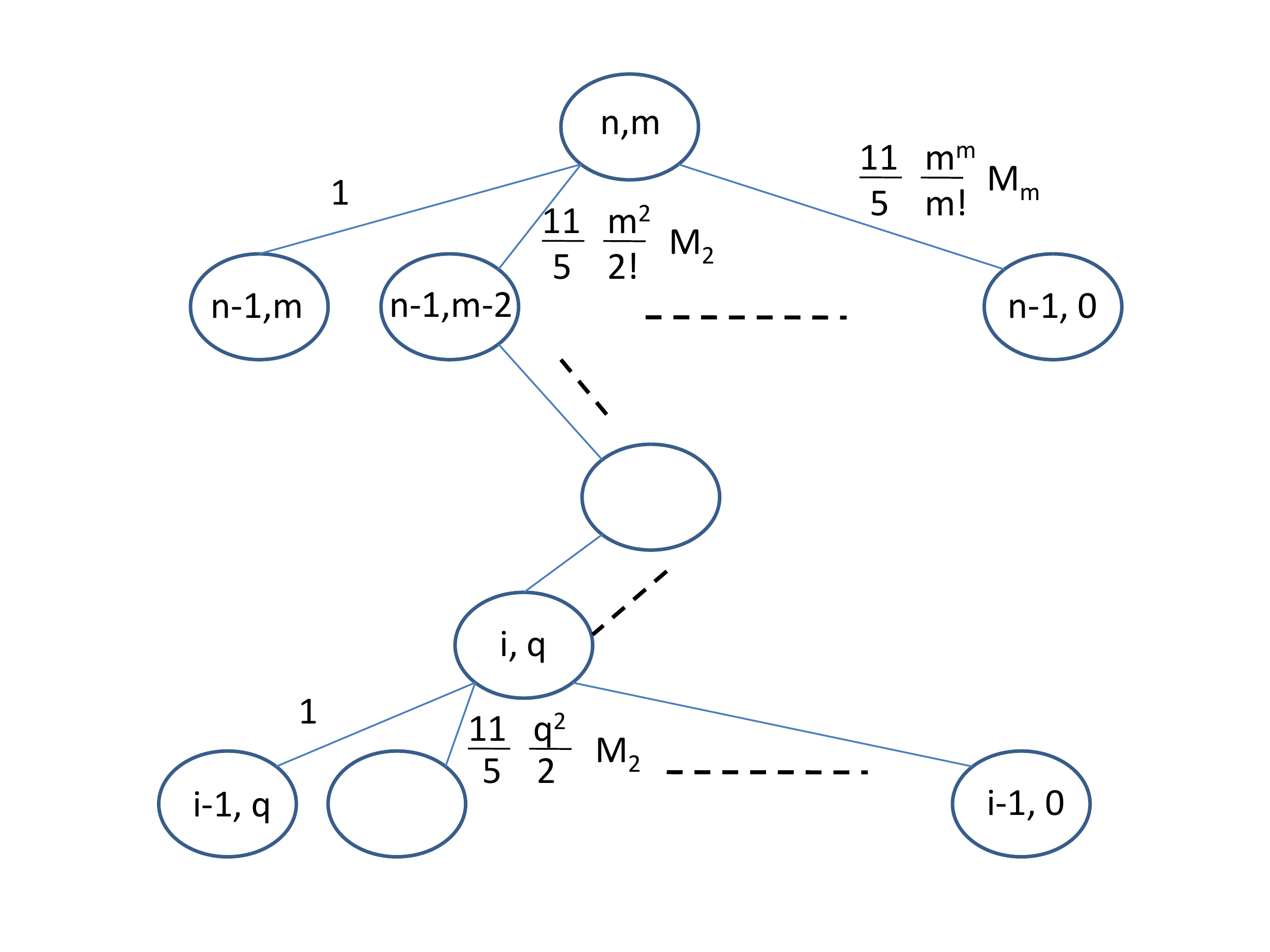}\\
\caption{Recursion Tree}\label{}
\end{center}
\end{figure}
constructed as follows:  The graph has a root marked $f(n,m)$. The root
has $(m/2)+1$ directed edges out of it going to $(m/2)+1$ nodes marked (respectively)
$f(n-1,m),f(n-1,m-2),\ldots f(n-1,0)$. The edges have weights associated with them which are
(respectively) $1,{11\over 5}{m^2\over 2!}M_2,{11\over 5}{m^4\over 4!}M_4,\ldots
{11\over 5}{m^m\over m!}M_m$. In general,
a node of the directed graph marked $f(i,q)$
(for $i\geq 2$, $0\leq q\leq m$, even) has $(q/2)+1$ edges going from it to nodes marked
$f(i-1,q),f(i-1,q-2),\ldots f(i-1,0)$; these edges have
``weights'' respectively $1,{11\over 5}{q^2\over 2!}M_2,{11\over 5}{q^4\over 4!}M_4,\ldots
{11\over 5}{q^q\over q!}M_q$ which are respectively at most
$$1,{11\over 5}{m^2\over 2!}M_2,{11\over 5}{m^4\over 4!}M_4,\ldots
{11\over 5}{m^q\over q!}M_q.$$ A node marked $f(1,q)$ has one child - a leaf marked
$f(0,0)$ connected by an edge of weight $M_q$. Define the weight of a path
from a node to a leaf as the product of the weights of the edges along the path.
It is easy to show by
induction on the depth of a node that $f(i,q)$ is the sum of weights of all
paths from node marked $f(i,q)$ to a leaf. [For example, if the assertion holds for all $i\leq n$,
then (\ref{evenM-special}) implies that it holds for the root.] We do not formally prove this here. A similar (slightly more complicated) Lemma - Lemma (\ref{recursive1})- will be proved during the proof of the Main Theorem.

Now, there is a 1-1 correspondence between paths from $f(n,m)$ to a leaf and elements of the
following set :
$L=\{ (l_1,l_2,\ldots l_n) : l_i\geq 0,\hbox{ even }; \sum_{i=1}^n l_i=m\}$; $l_i$ indicates that at level $i$
we take the $l_i$ th edge - i.e., we go from node $f(i,m-l_n-l_{n-1}-\ldots l_{i+1})$ to $f(i-1,m-l_n-l_{n-1}-\ldots l_i)$ on this path.
For an $l=(l_1,l_2,\ldots l_n)\in L$ and $t\in \{ 0,2,4,\ldots m\}$, define
$$g_t(l)=\hbox{ number of $i$ with $l_i=t$ }.$$
Clearly, the vector $g(l)=(g_0(l),g_2(l),\ldots , g_m(l))$ belongs to the set
$$H=\{ h=(h_0,h_2,h_4,\ldots h_m):\sum_{t} th_t=m;h_t\geq 0; \sum_t h_t=n\}.$$
Since the weight of an edge corresponding to $l_i$ at any level is at most
$({11\over 5})^zM_{l_i}{m^{l_i}\over l_i!}$, where $z=1$ iff $l_i\geq 2$, and the number of non-zero $l_i$ along any path is at most $m/2$, we have
\begin{align*}
f(n,m)&\leq \sum_{l\in L} ({11\over 5})^{m/2}\prod_t M_t^{g_t(l)} {m^{tg_t(l)}\over (t!)^{g_t(l)}}
\end{align*}
For an $h\in H$, the number of $l\in L$ with $g_t(l)=h_t\forall t$ is the number of ways of picking subsets of the $n$ variables of cardinalities $h_0,h_2,h_4,\ldots h_m$, namely,
$${n\choose h_0,h_2,h_4,\ldots h_m}={n!\over h_0!h_2!h_4!\ldots h_m!}\leq {n^{h_2+h_4+\ldots h_m}\over
h_2!h_4!\ldots h_m!}.$$
Thus, we have (using the assumed upper bound on conditional moments)
\begin{align}\label{7000}
&f(n,m)\leq ({11\over 5})^{m/2}\sum_{h\in H} {n^{h_2+h_4+\ldots h_m}\over
h_2!h_4!\ldots h_m!}\prod_t m^{th_t}{n^{h_t((t/2)-1)}\over m^{h_t((t/2)-1)}}\nonumber\\
&\leq ({11\over 5})^{m/2}\sum_h (nm)^{\sum_t th_t/2}{m^{h_2+h_4+\ldots +h_m}\over h_2!h_4!\ldots h_m!}\nonumber\\
&\leq
   ({11\over 5}nm)^{m/2}|H|\Max_{h\in H}{(em)^{h_2+h_4+\ldots +h_m}\over h_2^{h_2}h_4^{h_4}\ldots h_m^{h_m}} ,
   \end{align}
using Stirling inequality for factorial.
Now we will show that the maximum is attained when $h_2=m/2$ and the other $h_t$ are all zero.    In what follows $t$ only ranges over values $\geq 2$ for which $h_t\not= 0$.
\begin{align*}
&\prod_t {e^{h_t}m^{h_t}\over h_t^{h_t}}=\prod_t e^{h_t} \left( {m\over th_t}\right)^{h_t}t^{h_t}
\leq \prod e^{h_t}\left( 1+ \left( {m\over th_t}-1\right)\right)^{h_t}t^{h_t}
\leq \exp\left( \sum_t h_t+{m\over t}-h_t+h_t\ln t\right),
\end{align*}
using $1+x\leq e^x$ for all real $x$. Now, the function $\sum_t \left( {m\over t}+h_t\ln t\right)$ (considered as a function of the $h_t$) is linear and so its maximum over the simplex - $h\geq 0; \sum_t th_t=m$ - is attained at an extreme point. Hence $\sum_t \left( {m\over t}+h_t\ln t\right)\leq\Max_t \left( {m\over t}+{m\ln t\over t}\right).$ Now considered as a function of $t$, ${m\over t}+{m\ln t\over t}$ is decreasing, so the maximum of this over our range is at $t=2$. Thus, we have
\begin{equation}\label{7001}
\prod_t {e^{h_t}m^{h_t}\over h_t^{h_t}}\leq (2e)^{m/2}.
\end{equation}
Now,
we bound $|H|$: each element of $H$ corresponds to a unique ${m\over 2}$-vector $(h_2,2h_4,4h_8,\ldots )$ with coordinates summing
to $m/2$. Thus $|H|$ is at most the number of partitions of $m/2$ into $m/2$ parts which is ${m\choose m/2}\leq 2^m$. Plugging this and (\ref{7001}) into (\ref{7000}), we get the moment bound in the theorem.

The bound on $m$ th moment of $\sum_iX_i$ in the theorem will be used in a standard
fashion to get tail bounds. For any $t$, by Markov inequality, we get from the theorem
$\prob( |\sum_i X_i|\geq t)\leq {(24nm)^{m/2}\over t^m}$. The right hand side is minimized
at $m=t^2/(cn)$. So since the hypothesis of the theorem holds for this $m$, we get the claimed
tail bounds.

\proofend

The following Corollary is a strengthening of Chernoff Bounds.

\begin{corollary}\label{chernoff-cor}
Suppose $X_1,X_2,\ldots ,X_n$ are real valued random variables, $\sigma$ a positive
real and $m$ an even positive integer such that
\begin{align*}
&E(X_i^k|X_1+X_2+\ldots +X_{i-1})\leq\sigma^2\quad &\for k \hbox{  even  },\quad k\leq m\\
&EX_i(X_1+X_2+\ldots +X_{i-1})^k\leq 0\quad &\for k\hbox{  odd  },\quad k\leq m.
\end{align*}
Then, $E(\sum_{i=1}^n X_i)^m\leq (c_2nm\sigma^2)^{m/2}$ and $\sum_{i=1}^nX_i$ has $N(0,n\sigma^2)$ tails upto Min$(n\sigma^2,\sqrt{mn}\sigma)$.
\end{corollary}

\proofstart Let $t\in [0, \Min (n\sigma^2,\sqrt{mn}\sigma)]$.
We will apply the theorem with $m$ equal to the even integer nearest to $t^2/(c_1n\sigma^2)$ for a suitable $c_1>2$.
Since $t\leq n\sigma^2$, it is easy to see that $\sigma^2\leq \sigma^k (n/m)^{(k/2)-1}$ for any even $k$, so
the hypothesis of the theorem applies to the set of random variables - $(X_1/\sigma),(X_2/\sigma),\ldots (X_n/\sigma)$. So from the theorem, we get that
$$E(\sum_{i=1}^n X_i)^m\leq (c_2nm\sigma^2)^{m/2}$$ and so by Markov, we get
$$\prob (|\sum_{i=1}^n X_i|\geq t)\leq \left( {c_2nm\sigma^2\over t^2}\right)^{m/2}.$$
Now choose $c$ suitably so that ${c_2nm\sigma^2\over t^2}\leq {1\over 2}$ and we get the Corollary.

\endproof

\begin{remark}\label{chernoff}

The set-up for Chernoff bounds is: $X_1,X_2,\ldots X_n$ are i.i.d. Bernoulli random variables with $EX_i=\nu$. For any $t\leq n\nu$ Chernoff bounds assert:
$\prob \left( \left| \sum_{i=1}^n (X_i-\nu)\right| >t\right)\leq e^{-ct^2/(n\nu)}.$
We get this from the Corollary applied to $X_i-\nu$, since
$E(X_i-\nu)^2\leq\nu$ and
since $|X_i-\nu|\leq 1$, higher even moments of $X_i-\nu$ are at most the second moment. So, the hypothesis of the
Corollary hold with $\sigma^2=\nu$ and we can apply it.

The general Chernoff bounds deal with the case when the Bernoulli trials are independent, but not identical - $EX_i$ may be different for different $i$. This unfortunately is one of the points this simple theorem cannot deal with. However, the Main Theorem does deal with it and we can derive the general Chernoff bounds as a simple corollary of that theorem - see Remark (\ref{chernoff-2}).
\end{remark}

\section{Functions of independent random variables}\label{doob}

Theorem \ref{mainthm-special} and the Main Theorem (\ref{mainthm}) will often be applied to a real-valued
function $f(Y_1,Y_2,\ldots Y_n)$ of independent (not necessarily real-valued) random variables $Y_1,Y_2,\ldots $ to show concentration of $f$. This is usually done using the Doob's Martingale construction
which we recall in this section. While there is no new stuff in this section, we will introduce notation used throughout the paper.

Let $Y_1,Y_2,\ldots Y_n$ be independent random variables. Denote $Y=(Y_1,Y_2,\ldots Y_n)$.
Let $f(Y)$ be a real-valued function of $Y$.
One defines the classical Doob's Martingale:
$$X_i=E(f|Y_1,Y_2,\ldots Y_i)-E(f|Y_1,Y_2,\ldots Y_{i-1}).$$
It is a standard fact that the $X_i$ form a Martingale difference sequence and so
(\ref{snc}) is satisfied. We will use the short-hand $E^if$ to denote $E(f|Y_1,Y_2,\ldots Y_i)$, so $$X_i=E^if-E^{i-1}f.$$

Let $Y^{(i)}$ denote the $n-1$-tuple of random variables $Y_1,Y_2,\ldots Y_{i-1},Y_{i+1},\ldots Y_n$ and suppose $f(Y^{(i)})$ is also defined. Let
$$\Delta_i=f(Y)-f(Y^{(i)}).$$
\begin{equation}\label{XZ}
\hbox{Then,  } X_i =E^i\Delta_i-E_{Y_i}\left( E^i\Delta_i\right),
\end{equation}
since $Y^{(i)}$ does not involve $Y_i$.
$f,Y_i,X_i,\Delta_i$ will all be reserved for these quantities throughout the paper.
We use $c$ to denote a generic constant which can have different values.

\section{Random TSP with Inhomogeneous, heavy-tailed distributions}

One of the earliest problems to be studied under Probabilistic Analysis \cite{steele} is the concentration of the length $f$ of the shortest Hamilton cycle through a set of $n$ points picked uniformly independently at random from a unit square. Similarly, Karp's algorithm for the problem \cite{karp2} was one of the earliest polynomial time algorithms for the random variant of a problem which is NP-hard in the worst-case; see also \cite{steele2}.
It is known that $Ef\in \Theta(\sqrt n)$ and that
$f$ has $N(0,1)$ tails.
This was proved after many earlier steps by Rhee and Talagrand \cite{rhee-tsp} and Talagrand's inequality yielded a simpler proof of this. But Talagrand's method works only for independent points; under independence, the number of points in any sub-region of the unit square follows Poisson distribution which has exponentially falling tails.
Here, we will give a  simple self-contained proof of the concentration result
for more general distributions (of number of points in sub-regions) than the Poisson. Two important points of our more general distribution are
\begin{itemize}
\item Inhomogeneity (some areas of the unit square having greater probability than others)
is allowed.

\item heavier tails (for example with power-law distributions) than the Poisson are allowed.

\end{itemize}
We divide the unit square into $n$ small squares, each of side $1/\sqrt n$. We will generate at random a set $Y_i$ of points in the $i$ th small square, for $i=1,2,\ldots n$. We assume that the $|Y_i|$ are independent, but not necessarily identical random variables. Once the $|Y_i|$ are chosen, the actual sets $Y_i$ can be chosen in any (possibly dependent) manner (subject to the cardinalities being what was already chosen.) This thus allows for collusion where points in a small square can choose to bunch together or be spread out in any way.

\begin{theorem}\label{tsp}
Suppose there is a fixed
$c_1\in (0,1)$, an even positive integer $m\leq n$, and an $\epsilon>0$, such that for $1\leq i\leq n$ and $1\leq l\leq m/2$,
$$\prob ( |Y_i|=0 )\leq c_1\quad ;\quad E|Y_i|^l\leq (O(l))^{(2-\epsilon )l}.$$
Suppose $f=f(Y_1,Y_2,\ldots Y_n)$ is the length of the shortest Hamilton tour through $Y_1\cup Y_2\cup\ldots Y_n$.
Then, $f$ has $N(0,1)$ tails upto $\sqrt m$.
\end{theorem}

\begin{remark}
If each $Y_i$ is generated according to a Poisson of intensity 1
(=Area of small square times $n$), then $E|Y_i|^l\leq l^l$ and so the conditions of the theorem are satisfied for all $m$ (with room to spare).
\end{remark}
\begin{remark}
Note that if the hypothesis hold only upto a certain $m$, we get normal
tails upto $\sqrt m$. So for example $|Y_i|$ can have power law tails and
we still get a result, whereas the older results require exponential tails.
\end{remark}
\proofstart
Order the small squares in $\sqrt n$ layers - the first layer consists of all squares touching the bottom or left boundary; the second layer consists of all squares which are 1 square away from the bottom and left boundary etc. until the last layer is the top right square (order within each layer is arbitrary.)  Fix an $i$. Let $S_i$ be the $i$ th square.
Let $\tau = \tau (Y_{i+1},\ldots Y_n)$ be the minimum distance from a point of $S_i$ to a point in $Y_{i+1}\cup \ldots Y_n$ and $\tau_0=\Min (\tau, 2\sqrt 2).$
$\tau_0$ depends only on $Y_{i+1},\ldots Y_n$. (So, $E^i\tau_0=E\tau_0$.)
We wish to bound $\Delta_i=f(Y)-f(Y^{(i)})$
(see notation in section (\ref{doob})). For this, suppose we had a tour ${\cal T}$ through $Y^{(i)}$. We can break this tour at a point in $Y_{i+1}\cup Y_{i+2}\cup Y_n$
(if it is not empty) closest to $S_i$, detour to $S_i$, do a tour of $Y_i$ and then return to ${\cal T}$. If $Y_{i+1}\cup Y_{i+2}\cup Y_n$ is empty, we just break ${\cal T}$ at any point and do a detour through $Y_i$. So, we have
\begin{align*}
\Delta_i &\leq \tau_0+\hbox{ dist. from a point in $S_i$ to a point in $Y_i$ + length of tour thro' $Y_i$}+\tau_0\\
&\leq 2\tau_0+O(1/\sqrt n)+f(Y_i).
\end{align*}
Since $\Delta_i\geq 0$, we get using (\ref{XZ})for any even $l$:
\begin{align}\label{TSP-X}
&-E_{Y_i}(E^i\Delta_i)\leq X_i\leq E^i\Delta_i
\Longrightarrow |X_i|\leq 2E\tau_0+O(1/\sqrt n)+f(Y_i)\Longrightarrow E^{i-1}X_i^l\leq c^l (E\tau_0)^l+{c^l\over n^{l/2}}+{c^l\over n^{l/2}}E|Y_i|^{l/2}
\end{align}
where the last step uses the following well-known fact \cite{steele}.
\begin{claim}\label{tsp-det}
For any square $B$ of side $\alpha$ in the plane and any set of $s$ points in $B$, there is a Hamilton tour through the points of length at most $c\alpha \sqrt s$.
\end{claim}
First focus on $i\leq n-100\ln n$. We will see that we can get a good bound on $E\tau_0$ for these $i$.
For any $\lambda\in [0,5\sqrt{\ln n}/\sqrt n]$, there is a square region $T_\lambda$ of side $\lambda$ inside $S_{i+1}\ldots S_n$ (indeed, inside the later layers) which touches $S_i$. So,
$\prob (\tau \geq \sqrt 2\lambda )\leq \prob (T_\lambda \cap (Y_{i+1}\cup\ldots Y_n)=\emptyset)\leq e^{-cn\lambda^2}$ by the hypothesis that $\prob (|Y_j|=0)< c_1<1$.
This implies that
\begin{align*}
E\tau_0 &\leq \Prob \left( \tau \geq 5\sqrt{\ln n}/\sqrt n\right) (2\sqrt 2) + E(\tau|\tau \leq 5\sqrt{\ln n}/\sqrt n)\\
&\leq {c\over \sqrt n}+\left(  E(\tau^2|\tau \leq 5\sqrt{\ln n}/\sqrt n)\right)^{1/2}\leq {c\over\sqrt n}+\left( \int_0^\infty \lambda e^{-cn\lambda^2}\right)^{1/2}\leq {c\over\sqrt n}.
\end{align*}
Plugging this and the fact that
that $E|Y_i|^{l/2}\leq (O(l))^{(2-\epsilon )(l/2)}\leq (O(l))^l$ into (\ref{TSP-X}), we get
$E^{i-1}X_i^l\leq {l^l\over n^{l/2}}$.
We now apply theorem (\ref{mainthm-special}) to $c_6\sqrt n X_i$, for $i=1,2,\ldots n-100\ln n$
to get
\begin{equation}\label{n-100lnn}
E\left( \sum_{i=1}^{n-100\ln n}X_i\right)^m\leq (cm)^{m/2}.
\end{equation}
Now, we consider $i\geq n-100\ln n+1$. All of these squares are inside a square of side $\sqrt{\ln n}/\sqrt n$.
So, we have
$|\sum_{i=n-100\ln n+1}^n X_i|\leq 2\sqrt 2+{c\sqrt{\ln n}\sqrt{\sum_{i=n-100\ln n+1}^n|Y_i|}\over n^{1/2}}$.
Now using $E\left(\sum_{i=n-100\ln n+1}^n |Y_i|\right)^{m/2}\leq c(\ln n)^{m/2} m^{m-\epsilon m}$,
we get $E\left( \sum_{i=1}^nX_i\right)^m\leq (cm)^{m/2}$ which by the usual argument via Markov
inequality, yields the tail bounds asserted.
\proofend

\section{Minimum Weight Spanning tree}

This problem is tackled similarly to the TSP in the previous section. We will get the same
result as Talagrand's inequality is able to derive, the proof is more or less the same as our proof for the TSP, except that there is an added complication because adding points does not necessarily increase the weight of the minimum spanning tree. The standard example is when we already have the vertices of an equilateral triangle and add the center to it.

\begin{theorem}\label{mwst}
Under the same hypotheses and notation as in Theorem (\ref{tsp}),
suppose $f=f(Y_1,Y_2,\ldots Y_n)$ is the length of the minimum weight spanning tree on $Y_1\cup Y_2\cup\ldots Y_n$.
$f$ has $N(0,1)$ tails upto $\sqrt m$.
\end{theorem}

\proofstart  If we already have a MWST for $Y\setminus Y_i$, we can again connect the point in $Y_{i+1},\ldots Y_n$ closest to $S_i$ to $S_i$, then add on a MWST on $Y_i$ to get a spanning tree on $Y$. This implies again that
$\Delta_i\leq \tau_0+{c\sqrt{|Y_i|}\over\sqrt n}.$
But now, we could have $f(Y)<f(\hat Y)$. We show that
\begin{claim}\label{deletMWST}
$\Delta_i\geq -c_{10}\tau_0-{c\sqrt{|Y_i|}\over\sqrt n}.$
\end{claim}
\proofstart
We may assume that $Y_i\not=\emptyset$. Consider the MWST $T$ of $Y$.
We call an edge of the form $(x,y)\in T : x\in Y_i, y\in Y\setminus Y_i$, with $|x-y|\geq c_9/\sqrt n$, a long edge and an edge $(x,y) \in T: x\in Y_i, y\in  Y\setminus Y_i$, with $|x-y|< c_9/\sqrt n$
a short edge. It is well-known that the degree of each vertex in $T$ is $O(1)$ (we prove a more complicated result in the next para), so there are at most $6|Y_i|$ short edges; we remove all of them and add a MWST on the non-$Y_i$ ends of them. Since the edges are short, the non-$Y_i$ ends all lie in a square of side $O(1/\sqrt n)$, so a MWST on them is of length at most
$O(\sqrt {|Y_i|}/\sqrt n)$ by Claim (\ref{tsp-det}).

We claim that there are at most $O(1)$ long edges - indeed if $(x,y),(w,z)$ are any two long edges with $x,w\in Y_i$, we have $|y-z|\geq |x-y|-{\sqrt 2\over\sqrt n}$, since otherwise, $(T\setminus
(x,y))\cup (y,z)\cup (x,w)$ would contain a better spanning tree than $T$. Similarly,
$|y-z|\geq |w-z|-{\sqrt 2\over\sqrt n}$. Let $x_0$ be the center of square $S_i$. The above implies that in the triangle $x_0,y,z$, we have $|y-z|\geq |x_0-y|-{6\over\sqrt n},|x_0-z|-{6\over\sqrt n}$. But $|y-z|^2=|y-x_0|^2+|z-x_0|^2-2|y-x_0||z-x_0|\cos (y,x_0,z)$. Assume without loss of generality that $|y-x_0|\geq |z-x_0|$. If the angle $y,x_0,z$ were less than 10 degrees, then we would have $|y-z|^2\leq |y-x_0|^2+|z-x_0|^2-1.8|y-x_0||z-x_0|< (|y-x_0|-0.4|z-x_0|)^2$ a contradiction. So, we must have that the angle is at least 10 degrees which implies that there are at most 36 long edges.

Let $a$ be the point in $Y_{i+1},\ldots Y_n$ closest to $S_i$ if $Y_{i+1}\cup \ldots \cup Y_n$ is non-empty; otherwise, let $a$ be the point in $Y_1\cup Y_2\cup \ldots Y_{i-1}$ closest to $S_i$. We finally replace each long edge $(x,y),x\in Y_i$ by edge $(a,y)$. This clearly only costs us $O(\tau_0)$ extra, proving the claim.

Now the proof of the theorem is completed analogously to the TSP. \proofend

\section{Chromatic Number of inhomogeneous random graphs}

Martingale inequalities have been used in different (beautiful) ways on the chromatic
number $\chi$ of an (ordinary) random graph $G(n,p)$, where each edge is
chosen independently to be in with probability $p$
(see for example \cite{shamir},\cite{bollobas}, \cite{bollobas2},\cite{frieze}, \cite{luczak}, \cite{alon}, \cite{naor}).

Here we study chromatic number in a more general model.
An inhomogeneous random graph - denoted $G(n,P)$ - has vertex set $[n]$ and a $n\times n$ matrix $P=\{ p_{ij}\}$ where $p_{ij}$ is the probability that edge $(i,j)$ is in the graph. Edges are in/out independently. Let $$p={ \sum_{i,j}p_{ij}\over {n\choose 2}}$$ be the average edge probability. Let $\chi=\chi(G(n,P)$ be the chromatic number. Since each node can change the chromatic number by at most 1, it is trivial to see that $\prob (|\chi-E\chi|\geq t)\leq c_1e^{-c_2t^2/n}$ by H-A. Here we prove the
first non-trivial result, which is stronger than the trivial one when the graph is sparse, i.e., when $p\in o(1)$.
\begin{theorem}\label{chromatic}
$\chi$ of $G(n,P)$ has $N(0,n\ln n\sqrt p)$ tails upto $n\sqrt p$.
\end{theorem}
\begin{remark}
Given only $p$, note that $\chi$ could be as high as $\Omega(n\sqrt p)$ : for example, $p_{ij}$ could be $\Omega(1)$ for $i,j\in T$ for some $T$ with $|T|=O(n\sqrt p)$ and zero elsewhere.
\end{remark}
\proofstart
Let $p_i=\sum_j p_{ij}$ be the expected degree of $i$.
Let
$$S=\{ i : p_i\geq n\sqrt p\}.$$
$|S|\leq 2n\sqrt p$. Split the $n-|S|$ vertices of $[n]\setminus S$ into $k=(n-|S|)\sqrt p$ groups $G_1,G_2,\ldots G_k$ by picking for each vertex a group uniformly at random independent of other vertices. It follows by routine application of Chernoff bounds that with probability at least 1/2, we have : (i) for each $i$, the sum of $p_{ij}, j\in $ (same group as $i$) $\leq O(\ln n)$ and (ii) $|G_t|\in O(\ln n/\sqrt p)$ for all $t$. We choose any partition of $[n]\setminus S$ into $G_1,G_2,\ldots G_k$ satisfying (i) and (ii) at the outset and fix this partition.
Then we make the random choices to choose $G(n,P)$. We put the vertices of $S$ into singleton groups - $G_{k+1},\ldots G_{k+|S|}$.

Define $Y_i$ for $i=1,2,\ldots k+|S| $ as the set of edges (of $G(n,P)$) in
$G_i\times (G_1\cup G_2\cup \ldots G_{i-1})$.  We can define the Doob's Martingale
$X_i=E(\chi |Y_1,Y_2,\ldots Y_i)-E(\chi |Y_1,Y_2,\ldots Y_{i-1})$.
First consider $i=1,2,\ldots k$.
Define $\Delta_i$ as in section \ref{doob}.
Let $d_j$
be the degree of vertex $j$ in $G_i$
in the graph induced on $G_i$ alone. $\Delta_i$ is at most
$\max_{j\in G_i}d_j+1$, since we can always color $G_i$ with this many additional
colors. $d_j$ is the sum of independent Bernoulli random variables with $Ed_j=\sum_{l\in G_i}p_{jl}\leq O(\ln n)$. By Remark (\ref{chernoff-2}), we have that $E(d_j-Ed_j)^l\leq \Max
( (cl\ln n)^{l/2}, (cl)^l)$. Hence, $E^{i-1}(\Delta_i^l)\leq (cl)^l+(cl\ln n)^{l/2}$.

We will apply Theorem (\ref{mainthm-special}) to the sum
$${c_7X_1\over \ln n}+{c_7X_2\over\ln n}+\ldots {c_7X_k\over\ln n}.$$
It follows from the above that these satisfy the hypothesis of the Theorem provided $m\leq k$
. From this, we get
that
$$E\left(\sum_{i=1}^k X_i\right)^m\leq (cmk\ln n)^{m/2} .$$
For $i=k+1,\ldots k+|S|$, $\Delta_i$ are absolutely bounded by 1, so by the Theorem
$E( X_{k+1}+X_{k+2}+\ldots X_{k+|S|} )^m\leq (c|S|m)^{m/2}$. Thus,
$$E\left(\sum_{i=1}^{k+|S|} X_i\right)^m\leq (cmk\ln n)^{m/2}.$$

Let $t\in (0,n\sqrt p)$. We take
$m=$ the even integer nearest to $t^2/(c_4n\sqrt p \ln n)$ to get the theorem.
\proofend


\section{Random Projections}\label{JLsection}

A famous theorem of Johnson-Lindenstrauss \cite{vempala} asserts that if $v$ is picked uniformly at random from the surface of the unit ball in ${\bf R}^n$, then for $k\leq n$, and $\epsilon \in (0,1)$,
\footnote{A clearly equivalent statement talks about the length of the projection of a fixed unit length vector onto a random $k-$ dimensional sub-space.}
$\sum_{i=1}^k v_i^2$ has $N(0,\frac{k}{n^2})$ tails upto $\frac{k}{n}$.

The original proof
exploits the details of the uniform density
and simpler later proofs (\cite{arriaga}, \cite{dasgupta}, \cite{indyk})
use the Gaussian in the equivalent way of picking $v$.
Here, we will prove the same conclusion under weaker hypotheses which allows again longer tails
(and so does not use any special property of the uniform or the Gaussian). This is the first application which uses the Strong Negative Correlation condition rather than the Martingale Difference condition.

\begin{theorem}\label{jl}
Suppose $Y=(Y_1,Y_2,\ldots Y_n)$ is a random vector picked from a distribution such that (for a $k\leq n$)
(i) $E(Y_i^2|Y_1^2+Y_2^2+\ldots Y_{i-1}^2)$ is a non-increasing function of $Y_1^2+Y_2^2+\ldots Y_{i-1}^2$ for $i=1,2,\ldots k$ and (ii) for even $l\leq k$, $E(Y_i^l|Y_1^2+Y_2^2+\ldots Y_{i-1}^2)\leq (cl)^{l/2}/n^{l/2}$. Then, $\sum_{i=1}^k Y_i^2$ has $N(0,\frac{k}{n^2})$ tails upto $\frac{k}{n}$.
\end{theorem}

\proofstart
The theorem will be applied with $X_i=Y_i^2-EY_i^2$.
First, (i) implies for odd $l$:
$EX_i(X_1+X_2+\ldots X_{i-1})^l\leq 0$, by (an elementary version) of say, the FKG
inequality. [If $X_1+X_2+\ldots X_{i-1}=W$, then since $W^l$ is an increasing function of $W$
for odd $l$ and $E(X_i|W)$ a non-increasing function of $W$, we have $EX_iW^l=E_W\left( E(X_i|W)W^l\right)\leq E_W(E(X_i|W))EW^l=EX_iEW^l=0$.]
Now, for even $l$, $E^{i-1}(X_i^l) \leq 2^l EY_i^{2l}+2^l (EY_i^2)^l\leq (cl)^l/n^l$. So we may apply the theorem to the scaled variables $c_7nX_i$, for $i=1,2,\ldots k$ for $m\leq k$ to get that $n\sum_{i=1}^kX_i$ has $N(0,k)$ tails upto $O(\sqrt{kk})=O(k)$. So, $\sum_{i=1}^kX_i$ has $N(0,\frac{k}{n^2})$ tails upto $O(k/n^2)$ as claimed.
\proofend

{\bf Question} A common use of J-L is the following: suppose we have $N$ vectors $v_1,v_2,\ldots v_N$ is ${\bf R}^n$, where $n,N$ are high. We wish to project the $v_i$ to a space
of dimension $k<<n$ and still preserve all distances $|v_i-v_j|$. Clearly, J-L guarantees that for one $v_i-v_j$, if we pick a random $k$ dimensional space, its length is more or less preserved (within a scaling factor). Since the tail probabilities fall off exponentially in $k$, it suffices to take $k$ a polynomial in $\log N$ to ensure all distances are preserved. In this setting, it is useful to find more general choices of random subspaces (instead of picking them uniformly at random from all subspaces) and there has been some work on this (\cite{arriaga}, \cite{achlioptas}, \cite{ailon}). The question is whether Theorem 1 here or the Main Theorem can be used to derive more general results.

\section{Main Probability Inequality}\label{mainthcor}

Now, we come to the main theorem. We will again assume Strong Negative Correlation
(\ref{snc}) of the real-valued random variables $X_1,X_2,\ldots X_n$. The first main
point of departure from Theorem (\ref{mainthm-special}) is that we allow different variables
to have different bounds on conditional moments. A more important point will be that we
will use information on conditional moments conditioned on ``typical'' values of
previous variables as well as the pessimistic ``worst-case'' values.
More specifically, we assume
the following bounds on moments for $i=1,2,\ldots n$
($m$ again is an even positive integer):
\begin{equation}\label{Ml}
E(X_i^l|X_1+X_2+\ldots X_{i-1})\leq M_{il}\quad \for l=2,4,6,8\ldots m.
\end{equation}
In some cases, the bound $M_{il}$ may be very high
for the ``worst-case'' $X_1+X_2+\ldots X_{i-1}$. We
will exploit the fact that for a ``typical'' $X_1+X_2+\ldots X_{i-1}$,
$E( X_i^l|X_1+X_2+\ldots X_{i-1})$ may be much smaller. To this end, suppose
$${\cal E}_{i,l}\;\; , \; l=2,4,6,\ldots m\; ; i=1,2,\ldots n$$ are events.
${\cal E}_{i,l}$ is to represent the ``typical''
case. ${\cal E}_{1l}$ will be the whole sample space.
In addition to (\ref{Ml}), we assume that
\begin{align}
E( X_i^l|X_1+X_2+\ldots X_{i-1},{\cal E}_{i,l})&\leq L_{il}\label{Ll1}\\
\Prob ({\cal E}_{i,l})&=1-\delta_{i,l}\for l=2,4,6,8\ldots m\label{Li2}
\end{align}
Two quantities play a role in the theorem. The first is the ``average typical $l$ th moment'' $L_l$ which we define as
$$L_l=\frac{1}{n}\sum_{i=1}^n L_{i,l} \quad \for l=2,4,6,8\ldots m.$$
The second has to do with worst-case moments, but modulated by $\delta_{i,l}$. Let
$$\hat M_{i,l}=M_{i,l}\delta_{i,l}^{2/(m-l+2)}.$$
Note that while $M_{i,l}$ may be very large, one can make $\hat M_{i,l}$ smaller by controlling $\delta_{i,l}$.
\begin{theorem}[Main Theorem]\label{mainthm}
Let $X_1,X_2,\ldots X_n$ be real valued random variables satisfying Strong Negative
Correlation (\ref{snc}) and $m$ be a
positive even integer and $L_{l},\hat M_{il},\delta_{il}$ be as above. Then for $X=\sum_{i=1}^nX_i$,
\begin{align*}
&EX^m  \leq
(cm)^m\left(
\sum_{l=1}^{m/2}{1\over l^2}\left( {nL_{2l}\over m}\right)^{{1\over l}}\right)^{m/2}
+(cm)^{m}\sum_{l=1}^{m/2}
{1\over nl^2} \sum_{i=1}^n \left( n \hat M_{i,2l}\right)^{m/2l}.
\end{align*}
\end{theorem}
Besides the distinction between typical case and worst-case conditional moments which we already mentioned,
a second feature of the Theorem is similar to Theorem (\ref{mainthm-special}) in that the second moment term will often be the important one.
The $L$ term on the right hand side of the theorem is at most
$$(cm)^{m/2}\left( nL_2+\sqrt {nm}L_4^{1/2}+\ldots \right)^{m/2},$$
where we note that for $m<<n$, (which is the usual parameter setting with which the theorem will be applied) the coefficients of higher moments decline fast,
so that under reasonable conditions, the $nL_2$ term is what matters.
In this case, it will not be difficult to see that we get $N(0,nL_2)$ tails, as we would
in the ideal case when $X_i$ are independent and in the limit $X_1+X_2+\ldots +X_n$ behaves
like the normal (with variance equal to sum of the variances of the $X_i$, namely $nL_2$).

\begin{remark}\label{chernoff-2}
The general Chernoff bounds are a very special case: suppose $X_i, i=1,2,\ldots n$ are independent
Bernoulli trials with $EX_i=\nu_i$.
We will apply the theorem to bound the $m$ th moment of $X=\sum_i (X_i-\nu_i)$ and from that the tail probability. It is easy to see that $E(X_i-\nu_i)^l\leq \nu_i$ for all even $l$, so we may take $L_{i,2l}=\nu_i$ to satisfy the hypothesis of the Theorem {\bf for every $m$}. Let $\sum_i\nu_i=\nu$. We get $$EX^m\leq (cm)^{m/2}\left( m\sum_l (1/l^2){\nu^{1/l}\over m^{1/l}}\right)^{m/2}.$$ The maximum of $(\nu/m)^{1/l}$ occurs at $l=1$ if $\nu\geq m$ and at $l=m/2$ otherwise; in any case, it is at most $1+(\nu/m)$ and so we get (using $\sum_l (1/l^2)\leq 4$) for any $t>0$,
$$EX^m\leq (cm(\nu+m))^{m/2}\Longrightarrow \prob(|X|\geq t)\leq \left( { cm(\nu+m)\over t^2}\right)^{m/2}.$$
Now putting $m={t^2\over 2(\nu +t)}$, we get $\prob (|X|\geq t)\leq e^{-ct^2/(2(\nu+t))},$ which are Chernoff bounds.
\end{remark}

\section{Proof of the Main Theorem}

[The proof is complicated, not for lack of efforts on the part of the author. While certainly some of the intricate use of inequalities to get things to the final form which is usable may be necessary, it is possible that the reader may be luckier in simplifying the proof.]

We will use induction on $n,m$. At a general step of the argument, we will need to bound
$E\left( \sum_{i=1}^rX_i\right)^q$, where, $r\leq n$ and $q\leq m$, even. To bound this, let
$A=X_1+X_2+\ldots +X_{r-1}$. Binomial expansion gives us
$$E(\sum_{i=1}^rX_i)^q=E(A+X_r)^q = EA^q+qEX_rA^{q-1}+\sum_{l=2}^q {q\choose l}EX_r^lA^{q-l}.$$
The second term is non-positive by hypothesis. Also arguing exactly as in the proof of theorem (\ref{mainthm-special}), for odd $l\geq 3$,
$$EX_r^lA^{q-l}\leq {3\over 5} (EX_r^{l+1}A^{q-l-1}+EX_r^{l-1}A^{q-l+1}),$$
and so we get
\begin{equation}\label{thm2-recursion}
E\left(\sum_{i=1}^rX_i\right)^q\leq EA^q+3\sum_{{ l\geq 2\atop l\text{ even}}} {q\choose l}EX_r^lA^{q-l}.
\end{equation}
Without confusion, we will use ${\cal E}_{rl}$ to mean the 0-1 indicator variable of the event (defined earlier) ${\cal E}_{rl}$. Then, for even $l\geq 2$, we get
\begin{align*}
EX_r^lA^{q-l}&= EX_r^lA^{q-l}{\cal E}_{rl}+EX_r^lA^{q-l}(1-{\cal E}_{rl} )\leq L_{rl}EA^{q-l}+M_{rl}EA^{q-l}(1-{\cal E}_{rl})\\
&\leq L_{rl}EA^{q-l}+M_{rl}\left( EA^{q-l+2}\right)^{{q-l\over q-l+2}}\left( E(1-{\cal E}_{rl})\right)^{{2\over q-l+2}} \quad\text{  H\"older}\\
&\leq L_{rl}EA^{q-l}+M_{rl}\delta_{rl}^{{2\over m-l+2}}\left( EA^{q-l+2}\right)^{{q-l\over q-l+2}} \text{  since } m\geq q\\
&\leq L_{rl}EA^{q-l}+\left( \hat M_{rl}^{{2q\over l(q-l+2)}} (3m^2n^{2/l})^{{q-l\over q-l+2}}\right) \left( {\hat M_{rl}^{{(q-l)(l-2)\over l(q-l+2)}} \left( EA^{q-l+2}\right)^{{q-l\over q-l+2}}\over (3m^2n^{2/l})^{{(q-l)\over q-l+2}}}\right)\\
&\leq L_{rl}EA^{q-l}+\hat M_{rl}^{q/l}(3m^2n^{2/l})^{{q-l\over 2}}+{\hat M_{rl}^{{l-2\over l}}EA^{q-l+2}\over 3m^2n^{2/l}},
\end{align*}
where, in the last step, we have used Young's inequality which says that for any $a,b>0$ real and $s,r>0$ with ${1\over s}+{1\over r}
=1$, we have $ab\leq a^s+b^r$; we have applied this with $s=(q-l+2)/2$ and $r=(q-l+2)/(q-l)$.

Plugging this into (\ref{thm2-recursion}), we get:
$$E\left( \sum_{i=1}^r X_i\right)^q\leq \sum_{{l\geq 0\atop \hbox{even}}}^q a^{(q)}_{rl}E(\sum_{i=1}^{r-1}X_i)^{q-l},$$
\begin{align*}
a_{rl}^{(q)}&= 1+ 3{1\over 3m^2n}{q\choose 2}, & l=0\\
a_{rl}^{(q)}&= 3{q\choose l}L_{rl} +3{q\choose l+2}\hat M_{r,l+2}^{l/(l+2)}{1\over 3m^2n^{2/(l+2)}}  ,&  2\leq l\leq q-2\\
a_{rl}^{(q)}&= 3L_{rq}+3\sum_{{l_1\geq 2\atop \hbox{even}}}^q{q\choose l_1}\hat M_{rl_1}^{q/l_1} (3m^2)^{(q-l_1)/2}n^{(q-l_1)/l_1}  ,& l=q.
\end{align*}
It is easy to see that
\begin{align*}
a_{rl}^{(q)}\leq a_{rl}&= 1+ {1\over n}, \quad l=0\\
a_{rl}^{(q)}\leq a_{rl}&= 3{m\choose l}\left( L_{rl} + \hat M_{r,l+2}^{l/(l+2)}n^{-2/(l+2)}\right),\; 2\leq l\leq q-2\\
a_{rq}^{(q)}\leq \hat a_{rq}&= 3L_{rq}+3\sum_{{l_1\geq 2\atop \hbox{even}}}^q{q\choose l_1}\hat M_{rl_1}^{q/l_1}(3m^2)^{(q-l_1)/2}
                  n^{(q-l_1)/l_1}.\end{align*}
It is important to make $a_{r0}^{(q)}$ not be much greater than 1 because in this case only $n$ is reduced and so in the recurrence, this could happen $n$ times.
Note that except for $l=q$, the other $a_{rl}$ do not depend upon $q$; we have used $\hat a_{rq}$ to indicate that
this extra dependence. With this, we have
$$E\left( \sum_{i=1}^r X_i\right)^q\leq \hat a_{rq}+\sum_{{l\geq 0\atop \hbox{even}}}^{q-2} a_{rl}E(\sum_{i=1}^{r-1}X_i)^{q-l}.$$
We wish to solve these recurrences by induction on $r,q$.
Intuitively, we can imagine a directed graph with root marked $(n,m)$. The root has ${m\over 2}+1$
children which are marked $(n-1,m-l)$ for $l=0,2,\ldots m$; the node marked $(r,q)$ is trying to bound
$E(\sum_{i=1}^{r}X_i)^{q}$. There are also weights on the edges of $a_{rl}$. The graph keeps going
until we reach the leaves - which are marked $(1,*)$ or $(r,0)$. This is very similar to the recursion tree picture accompanying the proof of Theorem (\ref{mainthm-special}).
It is intuitively easy
to argue that the bound we are seeking at the root is the sum over all paths from the root to the leaves of the product
of the edge weights on the path. We formalize this in a lemma.

For doing that, for $1\leq r\leq n;2\leq q\leq m$, $q$ even and $1\leq i\leq r$ define
$S(r,q,i)$ as the set of $s=(s_i,s_{i+1},s_{i+2},\ldots s_r)$ with $s_i> 0;s_{i+1},s_{i+2},\ldots s_r\geq 0$ and $\sum_{j=i}^r s_j=q;s_j$ even.
\begin{lemma}\label{recursive1}
For any $1\leq r\leq n$ and any $q\leq m$ even, we have
$$E(\sum_{i=1}^rX_i)^q\leq \sum_{i=1}^r\sum_{s\in S(r,q,i)}\hat a_{i,s_i}\prod_{j=i+1}^r a_{j,s_j}.$$
\end{lemma}

{\bf Proof}
 Indeed, the statement is easy to prove for the base case of the induction - $r=1$ since
${\cal E}_{1l}$ is the whole sample space and $EX_1^q\leq L_{1q}$
. For the inductive step, we
proceed as follows.
\begin{align*}
&E(\sum_{i=1}^r X_i)^q\leq \sum_{{s_r\geq 0\atop\hbox{even}}}^{q-2} a_{r,s_r}E(\sum_{i=1}^{r-1}X_i)^{q-s_r}+\hat a_{r,q}\\
&\leq \hat a_{r,q}+
\sum_{i=1}^{r-1}\sum_{{s_r\geq 0\atop\hbox{even}}}^{q-2} a_{r,s_r}\sum_{s\in S(r-1,q-s_r,i)}\hat a_{i,s_i}\prod_{j=i+1}^{r-1} a_{j,s_j}.
\end{align*}
We clearly have
$S(m,q,m)=\{ q\}$ and for each fixed $i,1\leq i\leq r-1$, there is a 1-1 map\\
$S(r-1,q,i)\cup S(r-1,q-2,i)\cup \ldots S(r-1,2,i)\rightarrow S(r,q,i)$
given by\\
$ s=(s_i,s_{i+1},\ldots s_{r-1}) \rightarrow s'=(s_i,\ldots s_{r-1}, q-\sum_{j=i}^{r-1} s_j)$
and it is easy to see from this that we have the inductive step, finishing the proof of the Lemma.
\proofend

The ``sum of products'' form in the lemma
is not so convenient to work with. We will now get this to the
``sum of moments'' form stated in the Theorem. This will require a series of (mainly algebraic) manipulations with
ample use of Young's inequality, the inequality asserting $(a_1+a_2+\ldots a_r)^q\leq r^{q-1}(a_1^q+a_2^q+\ldots
a_r^q)$ for positive reals $a_1,a_2,\ldots $ and $q\geq 1$ and others.

So far, we have (moving the $l=0$ terms separately in the first step)
\begin{align}\label{sform}
&E\left( \sum_{i=1}^n X_i\right)^m
\leq \left( \prod_{i=1}^n a_{i0}\right) \sum_{i=1}^n \sum_{s\in S(n,m,i)} \hat a_{i,s_i}\prod_{{j=i+1\atop s_j\not= 0}}
   ^n a_{j,s_j}\nonumber\\
&\leq 3\sum_{i=1}^n \sum_{s\in S(n,m,i)} \hat a_{i,s_i}\prod_{{j=i+1\atop s_j\not= 0}}
   ^n a_{j,s_j}\nonumber\\
&\leq 3\sum_{t\geq 1}^{m/2}\left( \sum_{i=1}^n \hat a_{i,2t}\right) \sum_{s\in Q(m-2t)}\prod_{{j=1\atop s_j\not= 0}}^n a_{j,s_j}\\
&\hbox{     where,  } Q(q)=\{ s=(s_1,s_2,\ldots s_n) : s_i\geq 0\; \hbox{ even };\sum_j s_j=q\} \nonumber
\end{align}
Fix $q$ for now. For $s\in Q(q), l=0,1,2,\ldots q/2$, let $T_l(s)=\{ j: s_j=2l\}$ and $t_l(s)=|T_l(s)|$.
Note that $\sum_{l=0}^{q/2} lt_l(s)=q/2$.
Call $t(s)=(t_0(s),t_1(s),t_2(s),\ldots t_{q/2}(s))$ the ``signature''
of $s$. In the special case when $a_{il}$ is independent
of $i$, the signature clearly determines the ``$s$ term'' in the sum (\ref{sform}).
For the general case too, it will be useful
to group terms by their signature. Let (the set of possible signatures) be $T$. [$T$ consists of all $t=(t_0,t_1,t_2,\ldots t_{q/2})$ with
$t_l\geq 0\; \sum_{l=1}^{q/2} lt_l=q/2\; ;\sum_{l=0}^{q/2}t_l=n.$
\begin{align*}
\hbox{Now, }\sum_{s\in Q(q)}\prod_{{j=1\atop s_j\not= 0}}^n a_{j,sj}&= \sum_{t\in T}
    \sum_{T_0,T_1,T_2,\ldots T_{q/2}: |T_l|=t_l}
       ^{T_l\hbox{ partition }[n]}\prod_{l=1}^{q/2}\prod_{i\in T_l}a_{i,2l}  \\
&\leq \sum_{t\in T} \prod_{l=1}^{q/2}{1\over t_l!}
   \left( \sum_{i=1}^n a_{i,2l}\right) ^{t_l},
   \end{align*}
   since the expansion of $\left(\sum_{i=1}^n a_{i,2l}\right)^{t_l}$ contains $t_l!$ copies of
   $\prod_{i\in T_l}a_{i,2l}$ (as well other terms we do not need.)
Now define $ R=\{r=(r_1,r_2,\ldots r_{q/2}): r_l\geq 0;\sum_l r_l=q/2\}.$
We have
\begin{align}
&\sum_{t\in T} \prod_{l=1}^{q/2}
   {1\over t_l!}\left( \sum_{i=1}^n a_{i,2l}\right)^{t_l}\leq \sum_{r\in R}\prod_l {1\over (r_l/l)!}
   (\sum_{i=1}^n a_{i,2l} )^{r_l/l} \nonumber \\
&\leq {1\over (q/2)!}
\left( \sum_{l=1}^{q/2} m^{1-(1/l)}\left( \sum_ia_{i,2l}\right)^{1/l}
    \right)^{q/2},\label{tl}
\end{align}
where the first inequality is seen by substituting $r_l=t_ll$ and noting that the terms
corresponding to the $r$ such that $l|r_l\forall l$ are sufficient to cover the previous
expression and the other terms are non-negative. To see the second inequality, we just expand the
last expression and note that the expansion contains $\prod_l (\sum_i a_{i,2l})^{r_l/l}$ with
coefficient ${q/2\choose r_1,r_2,\ldots r_{q/2}}$ for each $r\in R$.
Now, it only remains to see that
$m^{r_l(1-(1/l))}\geq {r_l!\over (r_l/l)!}$, which is obvious.
Thus, we have plugging in (\ref{tl}) into
(\ref{sform}), (for some constant $c>0$; recall $c$ may stand for different constants at different points):
$$
EX^m\leq c^m
\sum_{t=1}^{{m\over 2}}{3\over ({m\over 2}-t)!}\left( \sum_{i=1}^n \hat a_{i,2t}\right)\left(
     \sum_{l=1}^{{m\over 2}-t} m^{1-{1\over l}} \left( \sum_i a_{i,2l}\right)^{{1\over l}}\right)^{{m\over 2}-t}. $$
\begin{align*}
\hbox{Now,  }({m\over 2}-t)!&\geq ({m\over 2}-t)^{{m\over 2}-t}e^{-{m\over 2}}e^t
\geq m^{{m\over 2}-t}e^{-{m\over 2}}
\Min_t \left[ \left(
  {{m\over 2}-t\over m}\right)^{{m\over 2}-t}e^t\right]\geq m^{{m\over 2}-t}(2e)^{-{m\over 2}},
\end{align*}
the last using Calculus to differentiate the log of the expression with respect to $t$ to see that the
min is at $t=0$. Thus, $$EX^m\leq c^m\sum_{t}\left[ \left({3\over m} \sum_{l=1}^{{m\over 2}-t} m^{1-{1\over l}} \left( \sum_i a_{i,2l}\right)^{{1\over l}}\right)^{{m\over 2}-t}\right]\left[ \sum_{i=1}^n \hat a_{i,2t}\right].$$
Let $\alpha,\beta$ denote the quantities in the 2 square brackets respectively. Young's inequality gives us:
: $\alpha\beta\leq \alpha^{m/(m-2t)}+\beta ^{m/2t}$. Thus,
\begin{equation}\label{aiaihat}
EX^m\leq \sum_{t=1}^{{m\over 2}}\left( \sum_i \hat a_{i,2t}\right)^{{m\over 2t}} +
    \sum_t \left( \sum_{l=1}^{{m\over 2}-1}m^{-{1\over l}}\left( \sum_i a_{i,2l}\right)^{{1\over l}}\right)^{{m\over 2}}
\end{equation}
In what follows, let $l_1$ run over even values to $m$ and $i$ run
from $1$ to $n$.
\begin{align}\label{aihat}
&\sum_{t=1}^{{m\over 2}}\left( \sum_i \hat a_{i,2t}\right)^{{m\over 2t}} \leq c^m
\sum_t \left( \sum_i L_{i,2t}\right)^{{m\over 2t}}\nonumber\\
&+c^mm^m\sum_t
\left( {1\over n}\sum_i \sum_{l_1\leq 2t}\left( {2t\over l_1m}\right)^{l_1} (n\hat M_{i,l_1})^{2t/l_1}\right)^{{m\over 2t}}\leq \nonumber\\
&c^m\sum_t (\sum_i L_{i,2t})^{{m\over 2t}}+ m^{m}\sum_{t,l_1}
{t^{{m\over 2t}} \over l_1^{l_1} n^{{m\over 2t}}}
\left( \sum_i (n\hat M_{i,l_1})^{2t/l_1}\right)^{{m\over 2t}}\nonumber\\
&\leq c^m\sum_t (\sum_i L_{i,2t})^{{m\over 2t}}+ c^mm^{m}\sum_{l_1}
{1\over nl_1^{l_1}}\sum_i (n\hat M_{i,l_1})^{m/l_1},
\end{align}
(using $t^{m/2t}\leq c^m$.)
\begin{align}\label{ai}
&\sum_{l=1}^{{m\over 2}-1}m^{-(1/l)}\left(\sum_i a_{i,2l}\right)^{{1\over l
}}\leq\nonumber\\
& \sum_{l=1}^{{m\over 2}-1}
 m^{-{1\over l}}\left( {m^{2l}\over (2l)!}\right)^{{1\over l}}
 \left( \sum_i L_{i,2l}+{\hat M_{i,2l+2}^{{l\over l+1}}\over n^{1/(l+1)}}
    \right)^{{1\over l}}\leq\nonumber\\
&m^2\sum_{l=1}^{{m\over 2}-1}
{m^{-{1\over l}}\over l^2}\left( \left( \sum_iL_{i,2l}\right)^{{1\over l}}
+\left(\sum_i \hat M_{i,2l+2}\right)^{{1\over l+1}}\right)\leq\nonumber\\
&m^2
\sum_{l=1}^{{m\over 2}}{m^{-{1\over l}}\over l^2}\left( \sum_i L_{i,2l}\right)^{1/l}+
m^2 \sum_{l=2}^{{m\over 2}}{1\over (l-1)^2}\left( \sum_i\hat M_{i,2l}\right)^{{1\over l}}.
\end{align}
We will further bound the last term using H\"older's inequality:
\begin{align}\label{lastterm}
&\left( \sum_{l=2}^{{m\over 2} }{\left( \sum_i\hat M_{i,2l}\right)^{1/l}\over (l-1)^2}\right)^{{m\over 2}}\leq\nonumber\\
&\left( \sum_{l=1}^\infty {1\over l^2}\right)^{(m-2)/2}\left(\sum_l {1\over (l-1)^2}\left(\sum_i \hat M_{i,2l}\right)^{{m\over 2l}}\right)\nonumber\\
&\leq 2^m\sum_{l=1}^{{m\over 2}}{1\over nl^2}\sum_i (n\hat M_{i,2l})^{{m\over 2l}}.
\end{align}
Now plugging (\ref{ai},\ref{aihat},\ref{lastterm}) into (\ref{aiaihat}) and noting that
$cm^{2-(1/l)}/l^2\geq 1$, we
get the Theorem.
\proofend

\section{Bin Packing}

Now we tackle bin packing. The input consists of $n$ i.i.d. items - $Y_1,Y_2,\ldots Y_n\in (0,1)$. Suppose $EY_1=\mu $ and $\var Y_1=\sigma$.
Let $f=f(Y_1,Y_2,\ldots Y_n)$ be the minimum number of capacity 1 bins into which the items $Y_1,Y_2,\ldots Y_n$ can be packed.
It was shown (after many successive developments) using non-trivial bin-packing theory (\cite{rhee}) that
$f$ has $N(0,n(\mu^2+\sigma^2))$ tails upto
$O(n(\mu^2+\sigma^2))$.
Talagrand \cite{talagrand} gives a simple proof of this from his inequality (this is the first of the six or so examples in his paper.) [We can also give a simple proof of this from our theorem.]

Talagrand \cite{talagrand} says (in our notation) ``especially when $\mu$ is small, one expects that the behavior of $f$ resembles the behavior of $\sum_{i=1}^nY_i$. Thereby, one should expect that $f$ should have tails of $N(0,n\sigma^2)$ or, at least, less ambitiously, $N(0,n(\mu^2+\sigma^2))$''.

However, $N(0,n\sigma^2)$ (as for sums of independent random variables) is easily seen to be impossible. An example is when items are of size $1/k$ or $(1/k)+\epsilon$ ($k$ a positive integer and
$\epsilon <<1/k$ is a positive real) with probability 1/2 each. $\sigma$ is $O(\epsilon)$.
It is clear that the number $n_1$ of $1/k$ items can be in ${n\over 2}\pm \Theta(\sqrt (n))$.
Now, a bin can have at most $k-1$ items if it has any $(1/k)+\epsilon$ item; it can have $k$ items if they are all $1/k$. Thus if $n_1$ number of $1/k$ items, we get
$$f = {n_1\over k}+{n-n_1\over k-1}+O(1) = {n\over 2}\left( {1\over k}+{1\over k-1}\right)\pm {\sqrt n\over k^2}.$$
From this it can be seen that the standard deviation of $f$ is $\Omega (\sqrt n\mu^2) >> \sqrt n\sigma$, establishing what we want.

Here we prove the best possible interval of concentration
when the items take on only one of a fixed finite set of values (discrete distributions - a case which has received much attention in the literature for example \cite{kenyon} and references therein).
[While our proof of the upper
bound here is only for problems with a fixed finite number of types, it would be
nice to extend this to continuous distributions.]

\begin{theorem}\label{binpack}
Suppose $Y_1,Y_2,\ldots Y_n$ are i.i.d. drawn from a discrete distribution with $r$ atoms, $r\in O(1)$, each with probability at least ${1\over \log n}$. Let $EY_1=\mu\leq
{1\over r^2\log n}$ and $\var Y_i=\sigma^2$. Then for any $t\in (0,n(\mu^3+\sigma^2))$, we have
$$\prob ( |f-Ef|\geq t+r)\leq c_1 e^{-c t^2/(n(\mu^3+\sigma^2))}.$$
Further, there is distribution for $Y_i$ in which $\var(f)\in\Omega(n(\mu^3+\sigma^2))$.
\end{theorem}
\proofstart
Let item sizes be $\zeta_1,\zeta_2,\ldots \zeta_j\ldots \zeta_r$ and the probability of picking type $j$ be $p_j$. [We will reserve $j$ to denote the $j$ th item size.] We have :
mean $\mu=\sum_j p_j\zeta_j$ and standard deviation
$\sigma=(\sum_j p_j (\zeta_j-\mu)^2)^{1/2}$.

Note that if $\mu\leq r/\sqrt n$, then earlier results already give concentration in an interval of length $O(\sqrt n(\mu +\sigma)$ which is then $O(r+\sigma)$, so there is nothing to prove. So assume that $\mu\geq r/\sqrt n$.

Define a ``bin Type'' as an $r-$ vector of non-negative integers
specifying number of items of each type which are together packable into one bin.
If bin type $i$ packs $a_{ij}$ items of type $j$ for $j=1,2,\ldots r$ we have
$\sum_j a_{ij}\zeta_j\leq 1$. Note that $s$, the number of bin types depends only
on $\zeta_j$, not on $n$.

For any set of given items, we may write a Linear Programming relaxation of the bin packing problem whose answers are within additive error $r$ of the integer solution. If there are $n_j$ items of size $\zeta_j$ in the set, the Linear program, which we call ``Primal'' (since later we will take its
dual) is :

Primal : ($x_i$ number of bins of type $i$.)
$$\Min \sum_{i=1}^s x_i\quad\hbox{ subject to } \sum_{i=1}^s x_ia_{ij}\geq n_j\forall j\; ; x_i\geq 0.$$
Since an optimal basic feasible solution has at most $r$ non-zero variables, we may just round these $r$ up to integers to get an integer solution; thus the additive error is at most $r$ as claimed.
In what follows, we prove concentration not for the integer program's value, but for the value of the Linear Program.
The Linear Program has the following dual :
$$\Max \sum_{j=1}^r n_jy_j\hbox{  s.t. }   \sum_j a_{ij}y_j\leq 1\; \for i=1,2,\ldots s;\; y_j\geq 0.$$
($y_j$ may be interpreted as the ``imputed'' size of item $j$)
Let $Y=(Y_1,Y_2,\ldots ,Y_n)$ and (for an $i$ we fix attention on) $Y'=(Y_1,Y_2,\ldots Y_{i-1},Y_{i+1},\ldots Y_n)$. We denote by $f(Y)$ the value of the Linear Program for the set of items $Y$. Let $\Delta_i=f(Y)-f(Y')$. The typical events ${\cal E}_i$ will just be that the number of copies of each $\zeta_j$
among $Y_1,Y_2,\ldots ,Y_{i-1}$ is close to its expectation:
$${\cal E}_i: \forall j, \left| \text{no.of copies of $\zeta_j$ in } Y_1,Y_2,\ldots ,Y_{i-1} - (i-1)p_j\right| \leq 100\sqrt{ m\ln {10m\over\mu}p_j(i-1)},$$
where $m$ is to be specified later, but will satisfy $\frac{1}{10}n(\mu^3+\sigma^2)$. We will use the Theorem with this parameter $m$.

We will make crucial use of the fact that second moments count highly for the bound in the theorem. So the main technical part of the proof is the following Lemma bounding typical conditional second moments.

\begin{lemma}\label{binpack-main}
Under ${\cal E}_i$, $\var (\Delta_i|Y_1,Y_2,\ldots ,Y_{i-1})\in O(\mu^3+\sigma^2)$.
\end{lemma}
\proofstart
Suppose now, we have already chosen all but $Y_i$.
Now, we pick $Y_i$ at random; say $Y_i=\zeta_k$. Let $Y=(Y_1,Y_2,\ldots Y_n)$ and
$Y'=(Y_1,Y_2,\ldots Y_{i-1},Y_{i+1},\ldots Y_n)$

Let
$$\Delta_i=f(Y)-f(Y').$$
Suppose we have the optimal solution of the LP for $Y'$.
There is a bin type which packs $\lfloor 1/\zeta_k\rfloor$
copies of item of type $k$; let $i_0$ be the index of this bin type. Clearly if we increase $x_{i_0}$ by
${1\over \lfloor 1/\zeta_k\rfloor}$, we get a feasible solution to the new
primal LP for $Y$. So
$0\leq \Delta_i\leq {1\over \lfloor 1/\zeta_k\rfloor}\leq
\zeta_k+2\zeta_k^2,$ which implies
\begin{align}
E(\Delta_i^2|Y')&\leq \sum_j p_j(\zeta_j+2\zeta_j^2)^2\leq \sum_jp_j\zeta_j^2+8\sum_jp_j\zeta_j^3\nonumber\\
&\leq\mu^2+\sigma^2+8\times 8 \sum_jp_j|\zeta_j-\mu|^3+8\times 8\sum_jp_j\mu^3 \leq \mu^2+65\sigma^2+64\mu^3\label{EZ2}.
\end{align}

Now, we lower bound $\Delta_i$ by looking at the dual. For this, let $y$ be the dual optimal
solution for $Y'$. (Note : Thus, $y=y(Y')$ is a function of
$Y'$.) $y$ is feasible to the new dual LP too (after adding in $Y_i$), since the dual constraints do not change.
So, we get:$\Delta_i\geq y_k$.
\begin{equation}
E(\Delta_i|Y')\geq \sum_j p_jy_j(Y')\label{111}.
\end{equation}
Also, recalling the bin type $i_0$ defined earlier, we see that
$y_k\leq 1/\lfloor (1/\zeta_k)\rfloor \leq \zeta_k+2\zeta_k^2$.
Say the number of items of
type $j$ in $Y'$ is $(n-1)p_j+\gamma_j.$
It is easy to see that $\zeta$ is a feasible dual solution. Since $y$ is an optimal
solution, we have
\begin{align*}
&\sum_j ((n-1)p_j+\gamma_j)y_j\geq \sum_j ((n-1)p_j+\gamma_j)\zeta_j.
\end{align*}
\begin{align}
&\mu-\sum_jp_jy_j=\sum_j p_j(\zeta_j-y_j)={1\over n-1}\left( \sum_j ((n-1)p_j+\gamma_j)(\zeta_j-y_j)\right) +{1\over n-1}\sum_j\gamma_j(y_j-\zeta_j)\nonumber\\
&\leq {1\over n-1}\sum_j\gamma_j(y_j-\zeta_j)\nonumber\\
&\leq {1\over n-1} (\sum_j (\gamma_j^2/p_j))^{1/2}(\sum_j p_j (y_j-\zeta_j)^2)^{1/2}\nonumber\\
&\leq {32(\mu+\sigma) r\over n}\Max_j |\gamma_j/\sqrt {p_j}|\label{delta},
\end{align}
where we have used the
fact that $-\zeta_j\leq y_j-\zeta_j\leq 2\zeta_j^2\leq 2\zeta_j$.
Let $(i-1)p_j+\gamma'_j$ and $(n-i)p_j+\gamma''_j$ respectively be the number of items of size $\zeta_j$ among $Y_1,Y_2,\ldots Y_{i-1}$ and $Y_{i+1},\ldots Y_n$. Since $\gamma_j''$ is the sum of $n-i$ i.i.d. random variables, each taking on value $-p_j$ with probability $1-p_j$ and $1-p_j$ with probability $p_j$, we have
$E(\gamma_j'')^2=\var (\gamma_j'')\leq np_j$.
Now, we wish to bound the conditional moment of $\gamma'_j$ conditioned on $Y_1,Y_2,\ldots Y_{i-1}$. But under the worst-conditioning, this can be very high. [For example, all fractions upto $i-1$ could be of the same type.] Here we exploit the typical case conditioning.  The expected number of ``successes'' in the $i-1$ Bernoulli trials is $p_j(i-1)$. By using Chernoff, we get (recall the definition of ${\cal E}_i$)
$\prob (\neg {\cal E}_i)= (\hbox{ say })\delta_i\leq \mu^{4m}m^{-4m}.$
Using (\ref{111}) and (\ref{delta}), we get
\begin{align*}
&E(\Delta_i|Y_1,Y_2,\ldots Y_{i-1};{\cal E}_i)\\
&\geq \mu-{32(\mu+\sigma) r\over n}E (\max_j {1\over\sqrt {p_j}}
            (100\sqrt{ m\ln (10m/\mu )p_j(i-1)} + (E(\gamma'')^2)^{1/2}))\\
&\geq \mu - c(\mu+\sigma)^{5/2}r\sqrt{\ln (10m/\mu)}-{c\mu r\over\sqrt n}\geq \mu - O(\mu^2),
\end{align*}
using $m\leq {1\over 10}n(\mu^3+\sigma^2)$.
So, we get recalling (\ref{EZ2}),
$$\var (\Delta_i|Y_1,Y_2,\ldots Y_{i-1};{\cal E}_i)=E(\Delta_i^2|Y_1,Y_2,\ldots Y_{i-1};{\cal E}_i) - (E(\Delta_i|Y_1,Y_2,\ldots Y_{i-1};{\cal E}_i))^2 \leq c(\mu^3+\sigma^2),$$
using ${r\over\sqrt n}\leq \mu\leq {1\over r^2\log n}$.
This completes the proof of the Lemma.
\proofend

Now, we have for the worst-case conditioning,
$$\var (\Delta_i|Y_1,Y_2,\ldots Y_{i-1})\leq E(\Delta_i^2|Y_1,Y_2,\ldots Y_{i-1})\leq c\mu^2.$$
We now appeal to (\ref{XZ}) to see that these also give upper bounds on $\var (X_i)$. As promised, dealing with higher moments is easy: note that $|\Delta_i|\leq 1$ implies that $L_{i,2l}\leq L_{i,2}$.
Now to apply the Theorem, we have
$L_{i,2l}\leq c(\mu^3+\sigma^2).$
So the ``$L$ terms'' are bounded as follows :
\begin{align*}
&\sum_{l=1}^{m/2}{m^{1-(1/l)}\over l^2}\left( \sum_{i=1}^n L_{i,2l}\right)^{1/l}
\leq \sum_{l=1}^{m/2} {m\over l^2}\left( {cn(\mu^3+\sigma^2)\over m}\right)^{1/l}\leq cn\left( \mu^3+\sigma^2\right)
\end{align*}
noting that $m\leq n(\mu^3+\sigma^2)$ implies that the
maximum of $((n/m)(\mu^3+\sigma^2 )^{1/l}$ is attained at $l=1$ and
also that $\sum_l (1/l^2)\leq 2$.
Now, we work on the $M$ terms in the Theorem.
$\max_{i}\delta_i\leq \mu^{4m}m^{-4m}=\delta^*$ (say).
$$\sum_{l=1}^{m/2}(1/n) \sum_{i=1}^n (n\hat M_{i,2l})^{m/2l}=
\sum_{l=1}^{m/2} e^{h(l)},$$
where $h(l) = {m\over 2l}\log n+{m\over l(m-2l+2)}\log \delta^*$.
We have $h'(l) = -{m\over 2l^2}\log n-\log\delta^* {m(m-4l+2)\over l^2(m-2l+2)^2}$.
Thus for $l\geq (m/4)+(1/2)$, $h'(l)\leq 0$ and so $h(l)$ is decreasing. Now for
$l < (m/4)+(1/2)$, we have
${m\over 2l^2}\log n \geq -(\log\delta^*){m(m-4l+2)\over l^2(m-2l+2)^2}$, so again $h'(l)\leq 0$. Thus, $h(l)$ attains its maximum at $l=1$, so
$
(36m)^{m+2}\sum_{l=1}^{m/2} e^{h(l)}\leq m(36m)^{m+3}n^{m/2}\delta^*$ giving us
$(36m)^{m+2}\sum_{l=1}^{m/2}(n\hat M^*_{2l})^{m/2l}
\leq (cnm(\mu^3+\sigma^2))^{m/2}.
$
Thus we get from the Main Theorem that
$E(f-Ef)^m\leq (cmn(\mu^3+\sigma^2))^{{m\over 2}},$
from which Theorem (\ref{binpack}) follows by the choice of $m=\lfloor {t^2\over c_5 n(\mu^3+\sigma^2)}\rfloor $.
\subsection{Lower Bound on Spread for Bin Packing}
This section proves the last statement in the theorem.Suppose
the distribution is :
\begin{align*}
\prob \left( Y_1={k-1\over k(k-2) }\right) &= {k-2\over k-1}\; ;\;
\prob \left( Y_1={1\over k}\right)= {1\over k-1}.
\end{align*}
This is a ``perfectly packable distribution'' (well-studied
class of special distributions) ($k-2$ of the large
items and 1 of the small one pack.) Also, $\sigma$ is small. But
we can have number of $1/k$ items equal to
${n\over k-1}-c\sqrt {{n\over k}}.$
Number of bins required $\geq \sum_i X_i = {n\over k}+{n\over
k(k-1)}+ c\sqrt {{n\over k}}\left({1\over k}\left( {k-1\over
k-2}-1\right)\right)\geq {n\over k-1}$. So at least
$c\sqrt{{n\over k}}$ bins contain only $(k-1)/k(k-2)$ sized items
(the big items). The gap in each such bin is at least $1/k$ for a
total gap of $\Omega (\sqrt n/k^{3/2})$. On the other hand, if the
number of small items is at least $n/(k-1)$, then each bin except
two is perfectly fillable.

\section{Longest Increasing Subsequence}

Let $Y_1,Y_2,\ldots Y_n$ be i.i.d., each distributed uniformly in $[0,1]$.
We consider here
$f(Y)=$ the length of the longest increasing subsequence (LIS)
of $Y$. This is a well-studied problem.
It is known that
$Ef=(2+o(1))\sqrt n$ (see for example \cite{aldous}). Since changing one $Y_i$ changes
$f$ by at most 1, traditional H-A yields $N(0,n)$ tails which is not so interesting.
Frieze \cite{frieze2} gave a clever argument (using a technique Steele \cite{steele} calls ``flipping'') to show concentration in intervals of length $n^{1/3}$. Talagrand \cite{talagrand} gave the first (very simple) proof of $N(0,\sqrt n)$ tails.
Here, we also supply a (fairly simple) proof from Theorem (\ref{mainthm}) of $N(0,\sqrt n)$ tails.
[But by now better intervals of concentration, namely $O(n^{1/6})$ are known, using detailed arguments specific to this problem \cite{bdj}.] Our argument follows from two claims below. Call $Y_i$ essential for $Y$ if $Y_i$ belongs to every LIS of $Y$ (equivalently, $f(Y\setminus Y_i)=f(Y)-1$.) Fix $Y_1,Y_2,\ldots Y_{i-1}$ and for $j\geq i$, let
$a_j=\prob\left( Y_j \hbox{ is essential for } Y|Y_1,Y_2,\ldots Y_{i-1}\right)$
\begin{claim}
$a_i,a_{i+1},\ldots a_n$ form a non-decreasing sequence.
\end{claim}
\proofstart
Let $j\geq i$. Consider
a point $\omega$ in the sample space where $Y_j$ is essential, but $Y_{j+1}$ is not. Map $\omega$ onto
$\omega'$ by swapping the values of $Y_j$ and $Y_{j+1}$; this is clearly a 1-1 measure preserving map. If
$\theta$ is a LIS of $\omega$ with $j\in\theta, j+1\notin \theta$, then $\theta \setminus j\cup j+1$ is an
increasing sequence in $\omega'$; so $f(\omega')\geq f(\omega)$. If $f(\omega')=f(\omega)+1$, then an LIS $\alpha$
of $\omega'$ must contain both $j$ and $j+1$ and so contains no $k$ such that $Y_k$ is between $Y_j,Y_{j+1}$.
Now $\alpha\setminus j$ is an LIS of $\omega$ contradicting the assumption that $j$ is essential for $\omega$. So
$f(\omega')=f(\omega)$. So, $j+1$ is essential for $\omega'$ and $j$ is not.
So, $a_j\leq a_{j+1}$.
\proofend
%
\begin{claim}
$a_i \leq c/\sqrt {n-i+1}$.
\end{claim}
\proofstart $a_i\leq {1\over n-i+1}\sum _{j\geq i}a_j$. Now $\sum_{j\geq i}a_j=a$ (say) is the expected
number of essential elements among $Y_i,\ldots Y_n$ which is clearly at most $Ef(Y_i,Y_{i+1},\ldots Y_n)\leq c\sqrt {n-i+1}$, so the claim follows.
\proofend

$\Delta_i$ is a 0-1 random variable
with $E(\Delta_i|Y_1,Y_2,\ldots Y_{i-1})\leq c/\sqrt {n-i+1}$. Thus it follows (using
(\ref{XZ}) of section (\ref{doob})) that
$$E(X_i^2|Y_1,Y_2,\ldots Y_{i-1})\leq c/\sqrt {n-i+1}.$$
Clearly, $E(X_i^l|Y_1,Y_2,\ldots Y_{i-1})\leq E(X_i^2|Y_1,Y_2,\ldots Y_{i-1})$
for $l\geq 2$, even. Thus we may apply the main Theorem with ${\cal E}_{il}$ equal to the
whole sample space. Assuming $p\leq \sqrt n$, we see that (using $\sum_l (1/l^2)=O(1)$)
$$E(f-Ef)^p\leq (c_1p)^{(p/2)+2}n^{p/4},$$
from which one can derive the asserted sub-Gaussian bounds.

\section{Number of Triangles in a random graph}

Let $f=f(G(n,p))$ be the
number of triangles in the random graph $G(n,p)$, where each edge is independently put in with probability $p$. There has been much work on the concentration of $f$.
\cite{kv1}, \cite{vu} discuss in detail why Talagrand's inequality cannot prove good
concentration when $p$ the edge probability is $o(1)$. [But we assume that
$np\geq 1$, so that $Ef=O(n^3p^3)$ is $\Omega(1)$.]
It is known (by a simple calculation - see \cite{jlr}
) that $$\var f=O(\Max (n^3p^3, n^4p^5)).$$
Our main result here is that $f$ has $N(0,\var f)$ tails upto $O^*((np)^{7/4})$, where, as usual, the $^*$ hides log factors. By a simple example, we see that $f$ does not have $N(0,\var f)$ tails beyond $(np)^{9/4}$. We note that our result is the first sub-Gaussian tail bound (with the correct variance) for the case when $p\leq 1/\sqrt n$. [For the easier case when $p=n^{-\alpha}, \alpha <1/2$, such a tail bound was known \cite{vu}, but only upto $(np)^{\epsilon}$ for a small $\epsilon >0$.]

The most
popular question about concentration of $f$ has been to prove upper bounds on
$\prob \left( f\geq (1+\epsilon) Ef\right)$ for essentially
$\epsilon \in \Omega(1)$ (see \cite{kv1}, \cite{janson}), i.e., for deviations as large as
$\Omega(Ef)$.  In a culmination of this line of work,
\cite{JOR} have proved that
$$\prob\left( f\geq (1+\epsilon) Ef\right) \leq ce^{-c\epsilon^2 n^{2}p^{2}}.$$
This is a special case of their
theorem on the number of copies of any fixed graph in $G_{n,p}$.
Their main focus is large deviations,
but for general $t$, putting $\epsilon=t/n^3p^3$ would only give us $e^{-t^2/(n^4p^4)}$. Also, \cite{kv1} develops a concentration inequality specially for polynomial functions of independent bounded random variables and \cite{vu} develops and surveys many applications of this inequalities; \cite{vu} discusses the concentration of the number of triangles as the ``principal example''.

\begin{theorem}\label{triangles}
$f$ has $N(0,\var f)$ tails upto $O^*((np)^{7/4})$.
\end{theorem}

\proofstart
Let $Y_i$ be the set of neighbors of vertex $i$ among $[i-1]$ and imagine adding
the $Y_i$ in order. [This is often called the vertex-exposure Martingale.]
We will also let $Y_{ij}$ be the 0-1 variable denoting whether there is an edge
between $i$ and $j$ for $j<i$.
The number of triangles $f$ can be written as $f=\sum_{i>j>k} Y_{ij}Y_{jk}Y_{ik}$.

As usual
consider the Doob Martingale difference sequence
$$X_i=E(f|Y_1,Y_2,\ldots Y_i)-E(f|Y_1,Y_2,\ldots Y_{i-1}).$$
It is easy to see that
$$X_i=\sum_{j<k\in [i-1]} Y_{jk}(Y_{ij}Y_{ik}-p^2)+(n-i)p^2\sum_{j<i}(Y_{ij}-p)=X_{i,1}+X_{i,2}\hbox{  (say)}.$$
Let $E^i$ denote $E(\cdot |Y_1,Y_2,\ldots Y_{i-1})$. We will be applying
our main concentration inequality Theorem (\ref{mainthm}) with $m=O(t^2/\var f)$.
Let $q$ be any even integer between 2 and $m$.
$E^i(X_i^q)\leq 2^qE^i(X_{i,1}^q)+2^qE^i(X_{i,2}^q)$.
Of the two, it is much easier to deal with $X_{i,2}$. Indeed we have using Corollary (\ref{chernoff-cor}).
\begin{equation}\label{Xi2}
E^i(X_{i,2}^q)\leq c^qn^qp^{2q}(npq)^{q/2}\leq (cn^3p^5q)^{q/2}.
\end{equation}
Let ${\cal E}_i$ be the event: (recall, as always, $c$ stands for poly$(\log n)$ and may have different values in different places)
\begin{align*}
{\cal E}_i:& |Y_j|\leq cnp\; \for j\leq i\\
&\forall S\subseteq [i-1], \text{ with }
|S|\leq cnp, \text{ we have } \sum_{j,k\in S} Y_{jk}\leq \max ( cn^2p^3,cnp)
\end{align*}
Now,
$$E^i(X_{i,1}^2)=\sum_{j_1<k_1<i}\sum_{j_2<k_2<i} Y_{j_1k_1}Y_{j_2k_2}E(Y_{ij_1}Y_{ik_1}-p^2)(Y_{ij_2}Y_{ik_2}-p^2).$$
Only terms where there are 2 or 3 distinct vertices among $j_1,j_2,k_1,k_2$ contribute to the expectation. The number of terms with 2 distinct vertices (and thus only one edge in $[i-1]$) is at most $n^2p$ under ${\cal E}_i$ and $E(Y_{ij_1}Y_{ik_1}-p^2)^2\leq p^2$, so the contribution of these terms is $O(n^2p^3)$. If there are 3 distinct vertices, we have a path of length 2 in $[i-1]$; there are $n^2p$ choices for the first edge of the path and $np$ choices of second edge under ${\cal E}_i$; finally, we have $|E (Y_{ij_1}Y_{ik_1}-p^2) (Y_{ij_1}Y_{ik_2}-p^2)| = O(p^3)$; so the total of these terms is $O(n^3p^5)$. Thus, we have
$$E^i(X_{i,1}^2|{\cal E}_i)\leq {c\var f\over n}.$$
Further, under ${\cal E}_i$, $|X_{i,1}|\leq a$, where, $a=\max (cn^2p^3,cnp)$ so we have for any even $l \geq 2$,
$E^iX_{i,1}^l\leq \var f a^{l-2}/ n$. We note also that $E^i(X_{i,2}^q)\leq
\var fa^{l-2}/n$, since $q\leq m\leq O^*(\sqrt {np})$ as is easy to see. Plugging these bounds into the ``$L$ terms'' of theorem (\ref{mainthm}), we get
$$\sum_{l=1}^{m/2} {1\over l^2} \left( {\sum_i L_{i,2l}\over m}\right)^{1/l}\leq ca^2
\sum_l {1\over l^2} \left( {\var f\over a^2m}\right)^{1/l}.$$
Since by the choice of $m$, we have $ma^2<\var f$, the maximum of $\left( {\var f\over a^2m}\right)^{1/l}$ is attained at $l=1$. Also $\sum_l (1/l^2)\in O(1)$. So, we have
\begin{equation}\label{L-terms}
m^m \left( \sum_{l=1}^{m/2} {1\over l^2} \left( {\sum_i L_{i,2l}\over m}\right)^{1/l}\right)^{m/2}\leq (m\var f)^{m/2}.
\end{equation}
Now, we bound the $M$ terms.
Since the expected number of edges within a particular $S\subseteq [i-1]$ with $|S|\leq cnp$ is
$O(n^2p^3)$, the probability that there are more than $\max (cn^2p^3,cnp)$
edges is most $e^{-cnp}$ for a particular $S$. Since there are at most $np {n\choose np}$
$S$ 's to consider, union bound gives us:
$$\delta_i=\prob (\neg{\cal E}_i)\leq e^{-cnp}.$$
We use a crude bound of $|X_i|\leq n^2$ to get $M_{i,l}\leq n^{2l}$. So,
$$(cm)^{m}\sum_{l=1}^{m/2}
{1\over nl^2} \sum_{i=1}^n \left( n \hat M_{i,2l}\right)^{m/2l}
\leq (cm)^m\sum_l (1/l^2)n^{m/2l}n^{2m}e^{-cnp/m}.$$
Again, it is easy to see that $m\leq O^*(\sqrt{np})$; so the above is at most
$(cm\var f)^{m/2}$. Together with the bound on the $L-$ terms, we now have
$$E(f-Ef)^m\leq (cm\var f)^{m/2},$$
from which the tail bound follows by using Markov as before.
\begin{remark}
It is easy to see that we do not have $N(0,\var f)$ tails beyond $(np)^{9/4}$ :
just take a random $G(n,p)$. Now add all the $(np)^{3/2}$ edges among the first
$(np)^{3/4}$ vertices; the probability of all these edges being present is
$e^{-c (np)^{3/2}}$ which is $e^{-t^2/n^3p^3}$, where the deviation $t$ from $Ef$
is $(np)^{9/4}$, namely the triangles among the first $(np)^{3/4}$ vertices.
\end{remark}
\begin{remark}
The inequalities in \cite{kv1} and \cite{vu} bound tails of polynomial functions of independent variables; the papers give many applications of them. Since most of the situations considered here are not polynomial functions, these are not applicable. But number of triangles is a polynomial of degree 3 in the underlying variables $Y_{ij}$ and so the main theorem of \cite{vu} (Theorem (4.2)) and Corollaries do apply. In that theorem, we have to choose $\tilde k=2$ or 3 and it is easy to see that with the conditions, we only get a tail bound which falls as $e^{-...t}$ and not $t^2$ as required for sub-Gaussian bounds.
\end{remark}

\section{Questions}

Many interesting open questions remain. Since the TSP is a classic problem, it would be interesting to strengthen/generalize results for the TSP. The first is to assume more limited independence: if one divides the unit square into $l$ pieces which have $Y_1,Y_2,\ldots ,Y_l$ as the set of points inside each respectively, can we prove concentration when $l\in o(n)$ and $E|Y_i|=n/l$ and assuming some moment conditions. Then, we have the question of extending concentration results under ``bursts in space'' to 3 and higher dimensions and finally, there are many other combinatorial problems \cite{steele} for which it would be interesting to prove such results.

We have not dealt much with ``bursts in time'', but the theorems here would seem to be applicable to such situations. In the bin-packing problem, it would be natural to assume that at each time $i$, one first picks the number of items which would arrive at that time and then have the items pick either adversarially or stochastically their sizes and prove concentration for the minimum number of bins. On-line versions of this problem are of interest. Queueing Theory has many examples of handling bursts and it remains to be seen how the results here may help in that area. 

The count of the number of not only triangles, but also other fixed graphs has been well-studied, but only for large deviations of the order of the expectation. It would be interesting to establish sub-Gaussian bounds as done here for triangles. This has some relation to the study of clustering coefficients and local communities in large (web-like) graphs.

\section{Comparisons with other inequalities}\label{comparisons}

The main purpose of this paper was to formulate and prove general probability inequalities which can be used to tackle the complicated combinatorial and other examples discussed. Here, we will compare our inequality to some others in the literature. For this we consider basic situations rather than complex ones to illustrate things better.

The ``sub-Gaussian'' behavior -
$e^{-{t^2.... }}$ with the ``correct'' variance (for example in Theorem (\ref{mainthm-special}) and Corollary
(\ref{chernoff-cor})) needs that the exponent of $m$ in the upper
bound in Theorem (\ref{mainthm-special}) be ${m\over 2}$.
Moment inequalities are of course well-studied and there are many sophisticated developments. One type of inequality is the Rosenthal type inequalities \cite{burkholder-2} which assert for Martingale difference sequence $X_1,X_2,\ldots X_n$ and even integer $m$:
$$E\left( \sum_{i=1}^n X_i\right)^m\leq f(m) \left( E \left( \sum_i E(X_i^2|X_1,X_2,\ldots X_{i-1})\right)^{m/2}  + E\max_i X_i^m\right).$$
Here, $f(m)$ has to be at least $cm/\ln m$ as shown by a simple example of
\cite{jsz}, which means that we cannot get sub-Gaussian bounds from these inequalities. The example is:
The $X_i$ are i.i.d. Bernoulli random variables with $X_i=1-(1/n)$ with probability $1/n$ and $-1/n$ with probability $1-(1/n)$ and $n=cm/\ln m$. For this, we have $E(\sum_{i=1}^nX_i)^m\geq n^m(1-(1/n))^m\prob (X_i=(1-(1/n))\forall i)\geq n^{m-o(m)}$. Our Theorem (\ref{mainthm-special}) can tackle the example: Note that for $l\geq \sqrt {\ln m}$, we have $(n/m)^{(l/2)-1} l!\geq 1$ and since $EX_i^l\leq 1/n$, the hypothesis of our Theorem (\ref{mainthm-special}) is satisfied. For $l\leq\sqrt{\ln m}$,
we see that $(n/m)^{(l/2)-1} l!\geq 1/n$ and this also suffices. So, our Theorem yields $E(\sum_i X_i)^m\leq (nm)^{m/2}$. But the example proves that $f(m)\geq (cm/\ln m)^{m}$.

Another class of inequalities are the Burkholder \cite{burkholder} type inequalities
which assert
$$E(X_1+X_2+\ldots +X_n)^m\leq g(m) E\left( X_1^2+X_2^2+\ldots +X_n^2\right)^{m/2},$$
for even integers $m$ when $X_i$ are Martingale differences. Here, since the right hand side involves taking the expectation of a power of the sum of $n$ quantities, we only gain if we could argue (in essence) that not many of them can be simultaneously high. Indeed, if we do not have any such information, then the best we might say is $X_1^2+X_2^2+\ldots +X_n^2\leq n\max_i|X_i|^2$, which only bounds the r.h.s. by $g(m)n^{m/2}E\max_iX_i^m$ and since it is known that $g(m)$ has to be at least $(cm)^{m/2}$, this does not give as strong results as Theorem (\ref{mainthm-special}). [The fact that $g(m)\geq (cm)^{m/2}$ follows from the simple example when $X_i$ are i.i.d., each equal to $\pm 1$ with probability 1/2 each.] But, here is a simple natural example where Burkholder inequality provably cannot derive something as strong as Theorem (\ref{mainthm-special}): let $Z_i$ be i.i.d., each Poisson with mean 1 and let $X_i =\pm Z_i, i=1,2,
\ldots n$, with probability 1/2 each, so $EX_i=0$. It is well-known that for even $l$, $EX_i^l=EZ_i^l=(cl)^l$, where $c$ here (and the rest of this section) involves constant and logarithmic (in $l$) factors. Theorem (\ref{mainthm-special}) directly yields $N(0,n)$ tails for $X=\sum_{i=1}^n X_i$ upto $n$. But to apply Burkholder, we must deal with $E(\sum_iX_i^2)^{m/2}$ for even $m$. We have
$$E\left(\sum_iX_i^2\right)^{m/2}\geq {n\choose m/2}(m/2)!(EX_1^2)^{m/2}+nEX_1^m\geq
  (cn)^{m/2}+(cm)^{m}.$$
So, the best one can ever prove is $EX^m\leq (cnm)^{m/2}+(cm)^{3m/2}$. Consider a tail probability $\prob (|X|\geq t)$; the best we could get for this from Burkholder type inequalities is
$$\prob (|X|\geq t)\leq {(cm)^{3m/2}\over t^m}+{(cnm)^{m/2}\over t^m}.$$
The minimum value of $(cm)^{3m/2}/t^m$ is easily seen by Calculus to be $e^{-ct^{2/3}}$ and when $n^{3/4}\in o(t)$, we have $t^{2/3}\in o(t^2/n)$, so we do not get $N(0,n)$ tails beyond $n^{3/4}$. One can ask if this is a cooked up example. But it occurs naturally - in many geometric probability results for example, where, $n$ i.i.d. points are picked uniformly from the unit square, it turns out that the ``Poisson approximation'' where instead one runs a Poisson process of intensity $n$ to get the points is more useful since, then, points in non-intersecting regions of the square are independent \cite{aldous1}. In this process, clearly the number of points in any region of area $1/n$ is Poisson with mean 1 and indeed, in our TSP and minimum weight spanning tree analysis, we used a generalization of this, allowing longer tails and dependence for the generation process and were still able to use Theorem (\ref{mainthm-special}).

The author has received many queries about how particular inequalities (the literature is clearly rich in this area with a number of clever papers, a majority appearing in the venerable journal: Annals of Probability) compares to the theorems here. An exhaustive comparison with each inequality in the literature would not be possible. But some more comparisons are given here. We consider three particular corollaries of our theorems - Generalized Chernoff bounds (GC) (Corollary (\ref{chernoff-cor}), Remarks (\ref{chernoff} and \ref{chernoff-2})), H-A and the Poisson example above. Our theorems can derive tail bounds for all of these.

A recent result on the line of Burkholder inequalities is for example, one in \cite{peligrad}, which asserts that
$$E(X_1+X_2+\ldots +X_n)^m\leq (cmn)^{m/2}EX_1^m,$$
where $m$ is again even and $X_i$ are now stationary Martingale differences. The $m^{m/2}$ is promising for getting sub-Gaussian bounds, but the high moment $EX_1^m$ on the right hand side means that Chernoff bounds don't follow from this. On the other hand, for stationary martingale differences, this is a strengthening of H-A. Talagrand's inequality can of course derive Chernoff bounds, but it only applies to independent random variables and so cannot derive H-A or GC. There are also inequalities based on the beautiful technique of Decoupling, for example Theorems 1.2A to 1.5B of \cite{lapena-2}. This works only for Martingale differences, requires bounds similar to our theorem (\ref{mainthm-special}), but for all moments, not just up to an $m$ precluding our Corollary (\ref{chernoff-cor}) and all other applications assuming only finite moments. But, we note that this does tackle the Poisson example and indeed, our theorem (\ref{mainthm-special}) is close in spirit to this, as discussed in Remark (\ref{decoupling}).
[Needless to add all the inequalities mentioned have their virtues which for want of space, we do not describe.] Here is a little table summarizing these comparisons.

\begin{tabular}{|c|c|c|c|}
  \hline
                 & \text{Poiss} & \text{H-A} & \text{GC} \\
   \text{THM 1}  & \tick     & \tick      & \tick \\
  \text{Rosent}  & X         & X          & X \\
  \text{Burkh}   & X         & \tick      & \tick \\
  \text{Decoup}  &  \tick    & \tick      & X \\
  \text{Talag}   & ???       & X          & X \\
  \hline
\end{tabular}

Legend: Poiss - the Poisson example above. GC - Generalized Chernoff.

Our crucial advantage is that while earlier moment inequalities generally do not focus on differentiating between the coefficients of different moments, the current paper
pays particular attention to the terms involving different moments. We are able to get a smaller coefficient on the higher moments which thus matter less; this is helpful, since lower moments are easier to bound tightly. This enables us to get the sub-Gaussian tails in the combinatorial situations discussed, whereas traditional inequalities do not get such bounds. It is worth noting that if we settle for an extra factor of $m^{m/2}$ in the bounds of our Theorem (\ref{mainthm}) (thus abandoning correct Gaussian tails) and also restrict only to Martingale differences instead of (\ref{snc}), then Burkholder's inequality would imply the theorem.

Another family of inequalities are the Efron-Stein type inequalities.
A recent result of Boucheron, Bousquet, Lugosi
and Massart \cite{lugosi} proves concentration for a real-valued
function $F$ of independent random variables
$Y_1,Y_2,\ldots Y_n$.
Let $Z=F(Y_1,Y_2,\ldots Y_n)$ and suppose
functions $Z_i=Z_i(Y_1,Y_2,\ldots Y_{i-1},Y_{i+1},\ldots Y_n)$ are arbitrary
functions. Their main
theorem is that
\begin{equation}\label{bblm}
E\left( (Z-EZ)_+\right)^m\leq (cm)^{m/2} E \left(\sum_{i=1}^n (Z-Z_i)^2\right)^{m/2}.
\end{equation}
[In the setting of independent random variables, this is in a way similar to Burkholder.]

Here, again, we sum up the
variations in $Z$ caused by all the $n$ variables and then take a high moment of it. The
advantage of this would be in situations where one can show that not too many of
individual $Y_i$
cause large changes for typical $Y_1,Y_2,\ldots Y_n$. [See \cite{lugosi}.] This
general line of
approach is also
reminiscent of Talagrand's inequality; but Talagrand allows
simultaneous change of variables. 
Note that (\ref{bblm})  has an exponent of $m/2$ on the $m$ which can lead
to the ideal sub-Gaussian behavior.
In contrast, our inequality (like Rosenthal's) only considers variations of one individual
variable at a time which
is in many cases easier to bound. We saw this in the case of Bin-Packing, coloring and
other examples. Even for the classical Longest Increasing Subsequence
(LIS) problem, where for
example, Talagrand's crucial argument is that only a small number $O(\sqrt n)$ of elements
(namely those in the current LIS) cause a decrease in the length of the LIS
by their deletion, we are able to bound individual variations (in essence
arguing that EACH variable has roughly only a $O(1/\sqrt n)$ probability
of changing the length of the LIS) sufficiently to
get a concentration result.

Note that if one can only handle individual variations,
then (\ref{bblm}) again essentially yields only
$$E\left( (Z-EZ)_+\right)^m\leq (cmn)^{m/2} \max_i E(Z-Z_i)^m.$$
In this case, arguments as in Theorem (\ref{mainthm-special})
as well as what we do for Bin-Packing and LIS which is
based mainly on the second moment, do not work, since the above
involves a high moment. There are many other specialized
ingenious probability inequalities in the literature; we have only
touched upon general ones.

Besides the situation like JL theorem, the Strong Negative correlation condition is also satisfied by the so-called  ``negatively associated'' random variables (\cite{joagdev},\cite{dubhashi}, \cite{boutsikas} for example).
Variables in occupancy (balls and bins) problems, 0-1 variables produced
by a randomized rounding algorithm of Srinivasan \cite{srinivasan} etc. are negatively associated.

{\bf Acknowledgements} Thanks to David Aldous, Alesandro Arlotto, Alan Frieze, Svante Janson, Manjunath Krishnapur, Claire Mathieu, Assaf Naor, Yuval Peres and Mike Steele, for helpful discussions.

\end{document}